\documentclass[12pt,a4paper]{amsart}
\usepackage{amssymb}
\usepackage{epsf}
\usepackage{latexsym}
\usepackage{pstricks}

\def\sym{\mathbb}
\def\maj{{\rm maj\,}}

\def\E{{\sym E}}
\def\N{{\sym N}}
\def\C{{\sym C}}

\def\K{{\sym K}}

\def\F{{\bf F}}
\def\G{{\bf G}}
\def\M{{\bf M}}
\def\S{{\bf S}}
\def\P{{\bf P}}
\def\NW{{\bf N}}
\def\WW{{\mathcal W}}

\def\ev{{\rm ev\,}}
\def\SG{{\mathfrak S}}

\def\pack{{\rm pack\,}}
\def\Des{{\rm Des\,}}
\def\Sym{{\bf Sym}}
\def\NCSF{{\bf Sym}}
\def\dim{{\rm dim\,}}
\def\Ker{{\rm Ker\,}}
\def\FQSym{{\bf FQSym}}
\def\WQSym{{\bf WQSym}}
\def\PQSym{{\bf PQSym}}
\def\CQSym{{\bf CQSym}}
\def\PBT{{\bf PBT}}
\def\DD{{\mathcal D}}
\def\btr{\blacktriangleright}

\newtheorem{example}{Example}[section]

\newtheorem{theorem}[example]{Theorem}

\newtheorem{corollary}[example]{Corollary}
\newtheorem{conjecture}[example]{Conjecture}

\newtheorem{proposition}[example]{Proposition}

\newtheorem{lemma}[example]{Lemma}

\def\up#1{\raise 1ex\hbox{\footnotesize#1}}

\def\Proof{\noindent \it Proof -- \rm}

\def\qed{\hspace{3.5mm} \hfill \vbox{\hrule height 3pt depth 2 pt width 2mm}}

\def\<{\langle}
\def\>{\rangle}

\def\ch{{\rm ch\,}}

\def\shuff#1#2{\mathbin{
      \hbox{\vbox{
        \hbox{\vrule
              \hskip#2
              \vrule height#1 width 0pt
               }%
        \hrule}%
             \vbox{
        \hbox{\vrule
              \hskip#2
              \vrule height#1 width 0pt
               \vrule }%
        \hrule}%
}}}

\def\shuffl{{\mathchoice{\shuff{7pt}{3.5pt}}%
                        {\shuff{6pt}{3pt}}%
                        {\shuff{4pt}{2pt}}%
                        {\shuff{3pt}{1.5pt}}}}%
\def\shuffle{\, \shuffl \,}

\def\gf#1#2{\genfrac{}{}{0pt}{}{#1}{#2}}

\def\vect{{\rm Vect\,}}
\def\rad{{\rm rad\,}}
\def\inv{{\rm inv\,}}

\title[The $(1-\E)$-transform]%
{The $(1-\E)$-transform in combinatorial Hopf algebras}

\author[F. Hivert, J.-G. Luque, J.-C. Novelli, J.-Y. Thibon]
{Florent Hivert, Jean-Gabriel Luque,\\ Jean-Christophe Novelli, and Jean-Yves
Thibon}

\address[Hivert and Luque]{LITIS, Universit\'e de Rouen ; Avenue de l'universit\'e ;
76801 Saint \'Etienne du Rouvray, France\\}

\address[Novelli and Thibon]{Institut Gaspard Monge, Universit\'e de
Marne-la-Vall\'ee \\
5, Boulevard Descartes \\Champs-sur-Marne \\77454 Marne-la-Vall\'ee cedex 2 \\
FRANCE}
\email[Florent Hivert]{hivert@univ-rouen.fr}
\email[Jean-Gabriel Luque]{luque@univ-rouen.fr}
\email[Jean-Christophe Novelli]{novelli@univ-mlv.fr}
\email[Jean-Yves Thibon]{jyt@univ-mlv.fr}

\date{\today}

\begin{document}
\begin{abstract}
We extend to several combinatorial Hopf algebras the endomorphism of symmetric
functions sending the first power-sum to zero and leaving the other ones
invariant. As a ``transformation of alphabets'', this is the
$(1-\E)$-transform, where $\E$ is the ``exponential alphabet'',
whose elementary symmetric functions are $e_n=\frac1{n!}$.
In the case of noncommutative symmetric functions, we recover
Schocker's idempotents for derangement numbers [Discr. Math. 269 (2003), 239].
From these idempotents, we construct subalgebras of the descent algebras
analogous to the peak algebras and study their representation theory. 
The case of $\WQSym$ leads to similar subalgebras of the Solomon-Tits algebras.
In $\FQSym$, the study of the transformation boils down to a simple solution
of the Tsetlin library in the uniform case.
\end{abstract}

\maketitle

\footnotesize
\tableofcontents
\normalsize

\section{Virtual alphabets in combinatorial Hopf algebras}

The notion of {\em virtual alphabet} provides a powerful symbolic notation
for dealing with endomorphisms of combinatorial Hopf algebras, at least
for those which can be realized in terms of (noncommutative) polynomials.
Examples include $\Sym$ (noncommutative symmetric
functions,~\cite{NCSF1,NCSF2}), $\FQSym$ (Free quasi-symmetric
functions~\cite{NCSF6,NCSF7}, which realize the Malvenuto-Reutenauer algebra
\cite{MR}), $\WQSym$ (Word quasi-symmetric functions~\cite{Hiv}), $\PQSym$
(Parking quasi-symmetric functions~\cite{NTp2}), and their subalgebras as long
as they are stable under the internal product.
These algebras can be regarded as generalizations of the Hopf algebra of
symmetric functions, for which this formalism was essentially promoted 
by A. Lascoux~\cite{LS,Las}. 

The algebra of symmetric functions is denoted by $Sym$. Apart from this
detail, our notation for symmetric functions follows Macdonald's
book~\cite{Mcd}.

\subsection{Virtual alphabets for symmetric functions}

It is convenient to use the generic term \emph{alphabet} to designate the
argument of a symmetric function. Indeed, a symmetric function $f$ is
characterized by its expression $f(x_1,x_2,\ldots)$ in terms of monomials in
an infinite sequence of independent indeterminates $x_i$, to which various
sets or multisets of algebraic expressions (numbers, monomials, cohomology
classes, vector bundles, \emph{etc}.) can be substituted.
The diversity of possible interpretations suggests to treat as far as possible
these arguments as formal symbols (or \emph{letters}, whence the term
\emph{alphabet}), and the possible occurence of multiple arguments suggests to
extend the usual meaning of this term and to understand it as a
\emph{multiset} of symbols.

Actually, a multiset $A=\{a,c,c,c,f,f\}$ is nothing but a formal linear
combination of symbols with nonnegative integer coefficients, and it is more
convenient to represent it as $A=a+3c+2f$. With this notation, the union of
multisets becomes just a sum $A+B$, and if we have sums, we can also have
differences $A-B$, at least when $B$ is contained in $A$.

This leads us to the point of this paragraph. When $B$ is not contained
in $A$, let us call the formal combination $A-B$ a \emph{virtual alphabet}.
It is easy to define the value of a symmetric function on a virtual alphabet.
Indeed, it is well known that any symmetric function $f$ can be expressed in
terms of the elementary symmetric functions $e_k$. The elementary functions
$e_k(A)$ of a genuine alphabet $A$ are the coefficients of the product
\begin{equation}
\lambda_t(A):=\prod_{a\in A}(1+ta)=\sum_k e_k(A)t^k
\end{equation}
and when $B$ is contained in $A$, those of the alphabet $C=A-B$ are
obtained by division of the generating functions:
\begin{equation}
\lambda_t(A-B)=\lambda_t(A)/\lambda_t(B)\,.
\end{equation}
When $B$ is not contained in $A$, this still defines coefficients
$e_k(A-B)$, which are, by definition, the elementary symmetric functions
of the virtual alphabet $A-B$. This defines all symmetric functions
of $A-B$, and one has for example the simple expression
$p_n(A-B)=p_n(A)-p_n(B)$ for the power sums.

More generally, a virtual alphabet can be defined by \emph{any} sequence
$(e_n(A))_{n\ge 1}$, interpreted as its elementary symmetric functions
(or as any sequence of independent generators of the algebra of
symmetric functions, such as power sums $p_n$ or complete homogeneous
functions $h_n$).

For example, although the exponential function $e^t$  has no zeros in the
complex plane, we can introduce a virtual alphabet $\E$ such that
\begin{equation}
e_n(\E) =\frac1{n!},\,\qquad\qquad \lambda_t(\E)=e^t\,.
\end{equation}
This is a useful trick, allowing to understand a lot of formulas in
combinatorics or analysis as specializations of simple identities on symmetric
functions. For example, the exponential generating function of the derangement
numbers $d_n$ is just
\begin{equation}
\sigma_t(1-\E):=\lambda_{-t}(1-\E)^{-1}=\frac{e^{-t}}{1-t}\,.
\end{equation}
Also, for a partition $\lambda$ of $n$, $n!s_\lambda(\E)=f_\lambda$, the
number of standard tableaux of shape $\lambda$.
Another virtual alphabet of interest is $1-q$, where $p_n(1-q)=1-q^n$. It
plays an essential role in the theory of Hall-Littlewood
functions~\cite{Mcd}.

Our expression of alphabets as formal sums or differences also
allows the consideration of products. For genuine alphabets $A$, $B$,
\begin{equation}
AB = \{a\,b\,\,|\,\, a\in A,\ b\in B\} =\sum_a a\sum_b b
\end{equation}
and for virtual alphabets, the symmetric functions of $AB$
are defined by any of the formulas like $p_n(AB)=p_n(A)p_n(B)$, or
\begin{equation}
h_n(AB)=\sum_{\lambda\vdash n}s_\lambda(A)s_\lambda(B)\,,
\end{equation}
the famous Cauchy identity (the $s_\lambda$ are the Schur functions).
More generally, for any function $f$, we have: $f(AX)=F(X)*\sigma_1(AX)$,
where $*$ denotes the internal product, and $\sigma_t(X)=\sum t^n h_n(X)$ is
the generating series of homogeneous complete functions.

Hence, with this formalism, we can consider \emph{transformations of
alphabets} on symmetric functions, which are ring homomorphisms mapping any
$f(X)$ to $f(\E X)$, $f((1-\E)X)$, $f((1-q)X)$, $f(X/(1-q))$, and so on.
Remark that in this setting, the virtual alphabet $1-q$ is the inverse of the
genuine alphabet $\{1,q,q^2\ldots\}\equiv 1+q+q^2+\cdots$.

Actually, $Sym$ is a $\lambda$-ring~\cite{Las,LS}, and the notion of plethysm
allows much more complicated transformations, but in the present paper, we
concentrate on the above mentioned ones since these can be defined only in
terms of the (combinatorial) Hopf algebra structure.

\subsection{Extension to combinatorial Hopf algebras}

Let $X$ and $Y$ be two infinite alphabets. Identifying expressions like
$f(X)g(Y)$ with the tensor product $f\otimes g$, we can regard the
transformations
\begin{equation}
\Delta:\ f(X)\longmapsto f(X+Y)
\end{equation}
and
\begin{equation}
\delta:\ f(X)\longmapsto f(XY)
\end{equation}
as linear maps $Sym\rightarrow Sym\otimes Sym$, that is, as comultiplications.
Moreover, both are (obviously) algebra morphisms, and the first one
being graded in the sense that
\begin{equation}
\Delta:\ Sym_n\longrightarrow \bigoplus_{i+j=n}Sym_i\otimes Sym_j
\end{equation}
defines actually a Hopf algebra structure, whose antipode is simply
$f(X)\mapsto f(-X)$.

Clearly, the transformations $f(\E X)$, $f((1-\E)X)$, $f((1-q)X$ and so on can
be defined only in terms of these coproducts, of the antipode, and of the
internal product. It turns out that most combinatorial Hopf algebras can be
realized in terms of polynomials in some infinite and \emph{totally ordered
alphabet}, denoted by $A=\{a_n|n\ge 1\}$ in the case of noncommuting letters,
and by $X=\{x_n|n\ge 1\}$ in the case of commuting letters. The basic example
is the pair Noncommutative Symmetric Functions -- Quasi-symmetric functions
$(\Sym,QSym)$ of mutually dual combinatorial Hopf algebras.  Sums and
differences of alphabets are defined on both sides, and it is possible to
make sense of the product $XA$ (see~\cite{NCSF2}).
Hence, our transformations are defined in this case.

The subject of this article is the study of the $(1-\E)$-transform in the
pair $(\Sym,QSym)$, and its extension to other combinatorial Hopf algebras.

\section{The $(1-\E)$-transform in $\Sym$}

\subsection{Background}

Our notations for the Hopf algebra of noncommutative symmetric functions are
as in~\cite{NCSF1,NCSF2}. This Hopf algebra is denoted by $\Sym$, or by
$\Sym(A)$ if we consider the realization in terms of an auxiliary alphabet.
Bases of its homogeneous component $\Sym_n$ are labelled by compositions
$I=(i_1,\ldots,i_r)$ of $n$. The noncommutative complete and elementary
functions are denoted by $S_n$ and $\Lambda_n$, and the notation $S^I$ means
$S_{i_1}\cdots S_{i_r}$. The ribbon basis is denoted by $R_I$.
The notation $I\vDash n$ means that $I$ is a composition of $n$.
The conjugate composition is denoted by $I^\sim$, the mirror image composition
by $\overline{I}$. The \emph{descent set} of $I$ is
$\Des(I) = \{ i_1,\ i_1+i_2, \ldots , i_1+\dots+i_{r-1}\}$.

The graded dual of $\Sym$ is $QSym$ (quasi-symmetric functions).
The dual basis of $(S^I)$ is $(M_I)$ (monomial), and that of $(R_I)$
is $(F_I)$.

The Hopf structures on $\Sym$ and $QSym$ allows one to partially extend the
$\lambda$-ring notation of ordinary symmetric functions (see~\cite{NCSF2},
and~\cite{Las} for background on the original commutative version).
If $A$ and $X$ represent totally ordered sets of noncommuting and commuting
variables respectively, the noncommutative symmetric functions of $XA$ are
defined by
\begin{equation}
\sigma_t(XA)=\sum_{n\ge 0}t^n S_n(XA) =
\prod_{x\in X}^{\rightarrow}\sigma_{tx}(A)= \sum_I t^{|I|}M_I(X)S^I(A)\,.
\end{equation}
Now, $X$ can be a virtual alphabet, defined by an arbitrary specialization of
an independ set of generators of $QSym$. 
An alternative way to express the transformation of alphabets defined by $X$
is~\cite{NCSF2}
\begin{equation}
F(XA)=F(A)*\sigma_1(XA)\,,
\end{equation}
where $*$ is the internal product. Since $\sigma_1(XA)$ is grouplike,
the $X$-transform is a bialgebra morphism, thanks to the splitting
formula
\begin{equation}
(F_1F_2\cdots F_r)*G = \mu_r \left[ (F_1\otimes\cdots\otimes F_r)* \Delta^r G
\right]
\end{equation}
where in the right-hand side, $\mu_r$ denotes the $r$-fold ordinary
multiplication and $*$ stands for the operation induced on $\Sym^{\otimes n}$
by $*$.

Thanks to the commutative image homomorphism $\Sym\rightarrow Sym$,
noncommutative symmetric functions can be evaluated on any element $x$ of a
$\lambda$-ring, $S_n(x)$ being $S^n(x)$, the $n$-th symmetric power. Recall
that $x$ is said \emph{of rank one} (resp.  \emph{binomial}) if
$\sigma_t(x)=(1-tx)^{-1}$ (resp.  $\sigma_t(x)=(1-t)^{-x}$).
The scalar $x=1$ is the only element having both properties. We usually
consider that our auxiliary variable $t$ is of rank one, so that
$\sigma_t(A)=\sigma_1(tA)$.

The argument $A$ of a noncommutative symmetric function can be a
\emph{virtual alphabet}. This means that, being algebraically independent,
the $S_n$ can be specialized to any sequence $\alpha_n\in{\mathcal A}$ of
elements of any associative algebra ${\mathcal A}$. Writing $\alpha_n=S_n(A)$
defines all the symmetric functions of $A$. Quasi-symmetric functions
of a virtual alphabet $X$ can be defined by a specialization of the
algebraic generators of $QSym$, which is more easily done by expressing
the noncommutative symmetric functions of $XA$ in terms of $A$.

The  specializations $X=\E$, defined by
\begin{equation}
S_n(\E A)=\frac{1}{n!}S_1(A)^n
\end{equation}
(so that $F_I(\E)=\frac1{n!}$) and $X=\frac{1}{1-q}$, for which
\begin{equation}
S_n\left(\frac{A}{1-q}\right)=\frac{1}{(q)_n}\sum_{I\vDash n}q^{\maj(I)}R_I(A)
\end{equation}
are of special importance. The second one can be used to define the
peak algebras and simplifies considerably their
investigation~\cite{BHT,ANT,NST}.
Here, we will study the $(1-\E)$- transform by the methods developed in these
references.

\subsection{Noncommutative derangement numbers}

A possible motivation for the  $(1-\E)$-transforms is the combinatorics of
derangements. Indeed, as already mentioned, the generating function of the
derangement numbers
\begin{equation}
D(t)=\sum_{n\ge 0} d_n\frac{t^n}{n!}=\frac{e^{-t}}{1-t} 
\end{equation}
can be expressed as
\begin{equation}
D(t)=\sigma_1((1-\E)t)\,.
\end{equation}
The specializations of the various bases of symmetric functions
at $1-\E$ are obtained by expanding the Cauchy kernel $\sigma_1((1-\E)X)$,
and, analogously, its quasi-symmetric functions can be defined as the
coefficients of the expansion of the non-commutative Cauchy kernel
\begin{equation}
\sigma_1((1-\E)A)=e^{-S_1(A)}\sigma_1(A)
\end{equation}
on any basis of noncommutative symmetric functions.

\bigskip
{\footnotesize
For example, since the quasi-monomial basis $M_I$ is dual to  $S^I$,
we have
\begin{equation} 
\sigma_1((1-\E)A)=\sum_{n\geq0} S_n((1-\E)A)= \sum_{I} M_I(1-\E) S^I(A),
\end{equation}
so that
\begin{equation} 
S_1((1-\E)A)= 0,\ S_2((1-\E)A)= S_2(A) - S^{11}(A)/2,\
\end{equation}
\begin{equation} 
S_3((1-\E)A)= S_3(A) - S^{21}(A) + S_{111}(A)/3,
\end{equation}
hence impliying
\begin{equation} 
M_{1}(1-\E)=0,\ M_{2}(1-\E)=1,\ M_{11}(1-\E)=-1/2,
\end{equation}
and
\begin{equation} 
M_{3}(1\!-\!\E)=1,\ M_{21}(1\!-\!\E)=-1,\
M_{12}(1\!-\!\E)=0,\ M_{111}(1\!-\!\E)=1/3.
\end{equation}
}

\bigskip

Since for $A=1$, $S_n((1-\E)A)= d_n/n!$, these noncommutative symmetric
functions might be called noncommutative derangement numbers
(see~\cite{NT3} for other examples of
\emph{noncommutative combinatorial numbers}).

\bigskip
{\footnotesize
Another natural noncommutative analogue of $d_n$ is given by
\begin{equation}\label{desar}
\lambda_{-t}(A)(1-tS_1(A))^{-1}
\end{equation}
which gives back the $d_n$ for $A=\E$, and $d_n(q)$ for $A=\frac1{1-q}$.
Its expansion of the ribbon basis is easily seen to be
\begin{equation}
\sum_{n\ge 0}t^n\sum_{i=1}^{[n/2]}\sum_{|J|=n-2i}R_{1^{2i}\triangleright J}
\end{equation}
(ribbons whose first column is of even height, or permutations whose
position of the first local minimum is even; the observation that these
permutations are counted by $d_n$ is due to D\'esarm\'enien, and an explicit
bijection has been given by D\'esarm\'enien and Wachs~\cite{DW1,DW2}).

Similarly,
\begin{equation}
\lambda_{u-t}(A)(1-tS_1(A))^{-1}
\end{equation}
gives a natural noncommutative analog of the generating series
$D(t,u)=e^{u-t}/(1-t)$ for permutations by number of fixed points,
and the expansion in $\FQSym$ of
\begin{equation}
\sigma_t(A) (1-t\,S_1(A))^{-1}  
\end{equation}
yields a class of permutations in bijection with arrangements: for
$A=\E$, this series is $\frac{e^t}{1-t}$.
}

\bigskip

Denote for short by a $\sharp$ the $(1-\E)$ transform, \emph{i.e.},
\begin{equation}
F^\sharp := F((1-\E)A)=F*\sigma_1^\sharp\,
\end{equation}
where
\begin{equation}
\label{eq-sigundiese}
\sigma_1^\sharp = \sigma_1((1-\E)A)=e^{-S_1}\sigma_1,
\end{equation}
that is,
\begin{equation}
\label{eq-sndiese}
S_n^\sharp = \sum_{i=0}^n (-1)^i \frac{S_1^i}{i!} S_{n-i}.
\end{equation}

\begin{lemma}
\label{lem-proj}
The $\sharp$-transform is a projector, \emph{i.e.},
\begin{equation}
\sigma_1^\sharp*\sigma_1^\sharp=\sigma_1^\sharp\,,
\end{equation}
or, equivalently,
\begin{equation}
S_n^\sharp * S_n^\sharp = S_n^\sharp \text{\rm \ \ for all $n$}.
\end{equation}
\end{lemma}

\Proof
Since the $\sharp$ transform is a homomorphism,
\begin{equation}
\begin{split}
\sigma_1^\sharp*\sigma_1^\sharp
&= (e^{-S_1}\sigma_1)*\sigma_1^\sharp 
= (e^{-S_1}*\sigma_1^\sharp)\, (\sigma_1*\sigma_1^\sharp) \\
&= e^{-S_1^\sharp}\, \sigma_1^\sharp = \sigma_1^\sharp.\\
\end{split}
\end{equation}
\qed

\bigskip
Hence, $S_n^\sharp$ corresponds to an idempotent $\delta_n$ of the descent
algebra $\Sigma_n$. We shall see later that the dimension of the
representation $\C\SG_n\delta_n$ is the derangement number $d_n$ and that
$\delta_n$ coincides with Schocker's derangement idempotent \cite{Sc}, whence
the name \emph{derangement algebra} below for the corresponding ideal of the
descent algebra.

\subsection{The small derangement algebra $\DD^{(0)}$}

\subsubsection{Definition of $\DD^{(0)}$}

Imitating the construction of the peak algebra in \cite{BHT}, we define the \emph{small
derangement algebra}, or the \emph{derangement ideal} as
\begin{equation}
\DD^{(0)} := \Sym^\sharp = \Sym(A)*\sigma_1^\sharp\,.
\end{equation}
By definition of the transformation of alphabet, $\DD^{(0)}$ is a Hopf
subalgebra of $\Sym$, and thanks to the second equality, each homogenous
component $\DD^{(0)}_n$ is a left ideal of $\Sym_n$ for the internal product,
explaing the two names for $\DD^{(0)}$, depending on which property
we want to emphasize.

\subsubsection{Dimension and bases of $\DD_n^{(0)}$}

We have $S_1^\sharp=0$, and the other $S_n^\sharp$ are clearly
algebraically independent. Hence, the dimension $d_n^{(0)}$ of $\DD_n^{(0)}$
is given by the number of compositions of $n$ with no part equal to $1$. These
elements satisfy the induction
\begin{equation}
d_0^{(0)}=1,\quad d_1^{(0)}=0,\quad
d_{n}^{(0)} = d_{n-1}^{(0)}+d_{n-2}^{(0)} \text{\ for $n\geq2$},
\end{equation}
hence are shifted Fibonacci numbers.

Let us say that a composition is \emph{non-unitary} if it does not contain the
part $1$. Then, the $S^{I\sharp}$ with $I$ non-unitary form a basis of
$\DD_n^{(0)}$.

Since the coproducts
\begin{equation}
\Delta S_n^\sharp=\sum_{i+j=n}S_i^\sharp\otimes S_j^\sharp
\end{equation}
are given by the same formula as those of the $S_n$, we see that
$\DD_n^{(0)}$ is isomorphic to the quotient $\Sym'$ of $\Sym$ by the 
two-sided ideal generated by $S_1$. Hence, its (graded) dual is isomorphic
to the subalgebra $QSym'$ of $QSym$ spanned by the $M_I$ for $I$ non-unitary.

Since $\sigma_1^\sharp$ is a projector, we have
\begin{equation}
f^\sharp*g^\sharp=(f*g)^\sharp\,.
\end{equation}
$\Sym'$ being stable under the internal product, $\DD_n^{(0)}$
is also $*$-isomorphic to $\Sym_n'$.

Recall that the ribbon basis of $\Sym$ is given by
\begin{equation}
R_I=\sum_{J\le I} S^J
\end{equation}
where $\le$ is the reverse refinement order. If $I$ is non-unitary,
so are all $J\le I$, thus
\begin{equation}
Q_I := R_I^\sharp\qquad (1\not \in I)
\end{equation}
is a basis of $\DD_n^{(0)}$.

\subsection{Algebraic structure of $(\DD_n^{(0)},*)$}

Since the construction of $\DD^{(0)}$ mimicks the one of the 
peak ideal, obtained in \cite{BHT} as the image of the
$(1-q)$-tranform at $q=-1$,
one may expect that $\DD^{(0)}$ shares many properties with the peak
ideal. There are some differences however.
Whilst the peak ideal has no unity for the internal product $*$, we have:
\begin{proposition}
\label{d0-unit}
For all $n$, $(\DD^{(0)}_n,*)$ is a unital algebra with $S_n^\sharp$
as neutral element.
\end{proposition}

\Proof
From Lemma \ref{lem-proj}, we already know that $\sigma_1^\sharp$ is neutral
on the right.
To prove that it is neutral on the left, let us consider its action
on the generating series of a basis of $\DD^{(0)}$
\begin{equation}
\sigma_1(X\cdot(1-\E)A)=\sum_I M_I(X){S^I}^\sharp(A)\,.
\end{equation}
We have
\begin{equation}
\sigma_1^\sharp*\sigma_1(X\cdot(1-\E)A)
=\sigma_1^\sharp*\sigma_1(XA)*\sigma_1^\sharp
\end{equation}
and
\begin{equation}
\begin{split}
\sigma_1^\sharp*\sigma_1(XA)
&= \left( e^{-S_1(A)}\sigma_1(A \right)) * \sigma_1(XA) \\
&= \mu\left[ (e^{-S_1(A)}\otimes \sigma_1(A))
       *_2 (\sigma_1(XA)\otimes\sigma_1(XA))\right] \\
&= (e^{-S_1(A)}*\sigma_1(XA))\,(\sigma_1(A) * \sigma_1(XA)) \\
&= e^{-S_1(XA)}\sigma_1(XA)
\end{split}
\end{equation}
so that
\begin{equation}
\begin{split}
\sigma_1^\sharp * \sigma_1(X\cdot(1-\E)A)
&= \left[ e^{-S_1(XA)}\sigma_1(XA) \right] * \sigma_1^\sharp \\
&= e^{-S_1(X\cdot(1-\E)A)}\,\sigma_1(X\cdot(1-\E)A) \\
&= \sigma_1(X\cdot(1-\E)A) \\
\end{split}
\end{equation}
since $S_1(X\cdot(1-\E)A)=0$. 
\qed

\subsection{Representation theory of $\DD_n^{(0)}$}

Now that we know that $\DD_n^{(0)}$ is unital, we can investigate its
representation theory. We first look at the idempotents.

\subsubsection{Idempotents in $\DD_n^{(0)}$}

Recall that the \emph{(right) Zassenhaus idempotents} $\zeta_n$ are
defined as the homogeneous elements of degree $n$ (that is,
$\zeta_n\in\NCSF_n$) satisfying (see~\cite{NCSF2}):
\begin{equation}
\sigma_1 =: \prod_{k\geq1} e^{\zeta_k}
         = e^{\zeta_1} e^{\zeta_2} e^{\zeta_3} \dots.
\end{equation}
Obviously, $\zeta_1=S_1$ so that $\zeta_1^\sharp=0$. Since the
$\sharp$-transform is multiplicative,
\begin{equation}
\sigma_1^\sharp= e^{\zeta_2^\sharp} e^{\zeta_3^\sharp} \dots
\end{equation}
but also
\begin{equation}
\sigma_1^\sharp= e^{-\zeta_1}\sigma_1=e^{\zeta_2} e^{\zeta_3} \dots
\end{equation}
so that
$\zeta_i^\sharp = \zeta_i$.
Extracting the term of degree $n$, we have
\begin{equation}
S_n^\sharp = \sum_{\gf{|\lambda|=n}{1\not\in\lambda}}
             \frac{\zeta^\lambda}{m_\lambda},
\end{equation}
where for a partition $\lambda=(\lambda_1\ge\lambda_2\ge\dots)$,
$\zeta^\lambda:=\zeta_{\lambda_1}\cdots \zeta_{\lambda_r}$,
and
$m_\lambda:=\prod_{i\geq1}m_i(\lambda)!$, where $m_i(\lambda)$ the
multiplicity of $i$ in $\lambda$.

We shall make use of the notation
$e_\lambda := \frac{\zeta^\lambda}{m_\lambda}$.
For a composition $I$, we set $\zeta^I=\zeta_{i_1}\dots\zeta_{i_r}$ and
$m_I=m_\lambda$ if $I\!\downarrow=\lambda$.

\subsubsection{A basis of idempotents}

The $\zeta_n$ are primitive elements with commutative image $p_n/n$,
hence are Lie idempotents~\cite{NCSF2}.
As with any sequence of Lie idempotents, we can construct an idempotent
basis of $\Sym_n$ from the $\zeta_n$. 

We need the following lemma from~\cite{NST}, easily derived from the splitting
formula (compare~\cite[Lemma 3.10]{NCSF2}).

Recall that the radical of $(\Sym_n,*)$ is
${\mathcal R}_n = {\mathcal R}\cap\Sym_n$, where ${\mathcal R}$ is the
kernel of the commutative image $\Sym\rightarrow Sym$. 

\begin{lemma}
\label{lem-NST}
Denote by $\SG(J)$ the set of distinct rearrangements of a composition $J$.
Let $I=(i_1,\ldots,i_r)$ and $J=(j_1,\ldots,j_s)$ be two compositions
of $n$. Then,

\smallskip (i) if $\ell(J)<\ell(I)$ then $\zeta^I*\zeta^J=0$.

\smallskip (ii) if $\ell(J)>\ell(I)$ then
$\zeta^I*\zeta^J \in \vect\<\zeta^K\, : \, K\in\SG (J)\>\cap{\mathcal R}$.
More precisely,
\begin{equation}
\zeta^I*\zeta^J 
= \sum_{\gf{\scriptstyle J_1,\ldots J_r}{\scriptstyle |J_k|=i_k}}
\< J\, ,\, J_1\shuffle \cdots \shuffle J_r \>\Gamma_{J_1}\cdots\Gamma_{J_r}
\end{equation}
where for a composition $K$ of $k$,
$\Gamma_K:=\zeta_k*\zeta^K$, and $\shuffle$ denotes the shuffle products
of compositions regarded as words over the positive integers.

\smallskip (iii) if $\ell(J)=\ell(I)$, then $\zeta^I*\zeta^J\not= 0$
only for $J\in\SG (I)$, in which case $\zeta^I*\zeta^J=m_I\, \zeta^I$. 
\end{lemma}
\hfill \qed

\begin{corollary}
\label{cor-NST}
(i) The elements
\begin{equation}
e_I=\frac1{m_I}\zeta^I\,,\quad I\vDash n\,,
\end{equation}
are all idempotents and form a basis of $\Sym_n$.
This basis contains in particular a complete set of minimal orthogonal
idempotents, $e_\lambda$ of $\Sym_n$.

(ii) The $e_I$ such that $I$ does not have a part equal to $1$ form a basis of
$\DD_n^{(0)}$.

(iii) The $e_\lambda$ with no part equal to $1$ in $\lambda$ form a complete
set of minimal orthogonal idempotents of $\DD^{(0)}_n$. \qed
\end{corollary} 

\subsection{Cartan invariants of $\DD^{(0)}_n$}

By (iii) of  lemma \ref{lem-NST}, the indecomposable projective module
$P_\lambda=\DD^{(0)}_n*e_\lambda$ contains the $e_I$ for $I\in\SG(\lambda)$.
For $I\not\in\SG(\lambda)$, (i) and (ii) imply that $e_I*e_\lambda$ is in
$\vect\<\zeta^K\, : \, K\in\SG (\lambda)\>$. Hence, this space coincides with
$P_\lambda$. So, we get immediately an explicit decomposition
\begin{equation}
\DD^{(0)}_n=\bigoplus_{\lambda\vdash n,\ 1\not\in\lambda}P_\lambda\,,\qquad
P_\lambda=\bigoplus_{I\in\SG(\lambda)}\C e_I\,.
\end{equation}
The Cartan invariants
\begin{equation}
c_{\lambda,\mu}=\dim\,( e_\mu*\DD^{(0)}_n*e_\lambda)
\end{equation}
are now easily obtained. The above space is spanned by the
\begin{equation}
e_\mu*e_I* e_\lambda = e_\mu* e_I\,,\quad I\in\SG(\lambda)\,.
\end{equation}

\noindent
From (ii) of Lemma~\ref{lem-NST}, this space has the dimension of the space 
$[S^\mu(L)]_\lambda$, spanned by all symmetrized products of Lie polynomials
of degrees $\mu_1,\mu_2,\ldots$ formed from $\zeta_{i_1},\zeta_{i_2},\ldots$,
as in the classical result of Garsia-Reutenauer for the descent
algebra~\cite{GaR}.

\bigskip
In the following examples, partitions are ordered by reverse lexicographic
order.
For $n\leq4$, the Cartan matrix of $\DD_n^{0}$ is trivial: it is the identity
matrix of size $d_n^{(0)}$,  since there is at most one value in a partition
of $n$ with no part one.
For $n$ up to $7$,  the Cartan invariants are given by the matrices
(\ref{qCar-dn0-567}) and~(\ref{qCar-dn0-89}) at $q=1$. Indeed, the
$q$-analogues defined from the Loewy series can be explicitely calculated.

\subsection{Quiver and $q$-Cartan invariants (Loewy series)}

We shall use the following modified refinement order on partitions: we write
$\lambda\prec_{p}\mu$, and say that $\lambda$ is $p$-finer than $\mu$, if each
part of $\mu$ is either a part of $\lambda$ or  a sum of distinct parts of
$\lambda$.  Hence, $\lambda$ covers $\mu$ iff $\mu$ is obtained from $\lambda$
by merging two distinct parts.

Still relying upon point (ii) of the lemma, we see that $c_{\lambda,\mu}=0$ if
$\lambda$ is not $p$-finer than (or equal to) $\mu$, and that if $\mu$ is
obtained from $\lambda$ by adding up two parts $\lambda_i,\lambda_j$, $e_\mu*
e_I=0$ if $\lambda_i=\lambda_j$ and is a nonzero element of the radical
otherwise.

The $q$-analogues of the Cartan invariants
\begin{equation}
c_{\lambda,\mu}
  := \sum_{k} q^k
     \dim\left[ (e_\mu*\rad^k\DD_n^{(0)}*e_\lambda)
                 / \rad^{k+1}\DD_n^{(0)} \right]
\end{equation}
can now be obtained from Proposition~\ref{d0-unit} and the following lemma.

\begin{lemma}
\label{BAe}
Let $A$ be an associative algebra. If $e$ is an idempotent of $A$ such that
$e$ is neutral in $B=Ae$, then
\begin{equation}
\rad^k B = e (\rad^k A) e.
\end{equation}
\end{lemma}

\Proof
Let $x\in \rad B$. There exists an integer $n$ such that
$(xB)^n=0$. Then,
\begin{equation}
(x eAe)^n = (xA)^n e = 0,
\end{equation}
so that, as well
\begin{equation}
(xA)^n exA = (xA)^n xA = (xA)^{n+1} =0.
\end{equation}
Thus, the right ideal $xA$ is nilpotent, which proves that $x\in\rad A$.
Since $x=exe$, $x\in e\,\rad A e$ and we have shown that
$\rad B\subseteq e\,\rad A e$.

Conversely, if $x\in\rad A$, so that $x^n=0$ for a certain $n$, then
$(exe)^n=x^ne=0$, whence $exe\in\rad B$, which proves the claim for $k=1$.

Now, if $x\in\rad^k B$, $x=x_1\dots x_k$ with $x_i=ey_ie$, for some
$y_i\in\rad A$. Hence
\begin{equation}
x = ey_1e\dots ey_ke = e y_1y_2\dots y_r e \in e\,\rad^k Ae.
\end{equation}
Conversely, any $x$ of the form $ey_1\dots y_re$ with $y_i\in\rad A$ can be
rewritten as $x=ey_1e\dots ey_ke\in \rad^k B$.
\qed

\medskip
Applying this to $A=\NCSF_n$, and $e=S_n^\sharp$, we obtain from the known
description of the $q$-Cartan matrices of $\NCSF_n$~\cite{Sal}:

\begin{theorem}
(i) In the quiver of $\DD^{(0)}_n$, there is an arrow $\lambda\rightarrow\mu$
iff $\mu$ is obtained from $\lambda$ by adding two distinct parts,

\noindent
(ii) The $q$-Cartan invariants of $\DD^{(0)}_n$ are given by
\begin{equation}
c_{\lambda,\mu}(q)= c_{\lambda,\mu}\, q^{\ell(\lambda)-\ell(\mu)}.
\end{equation}
if $\lambda$ is finer than or equal to $\mu$, and $c_{\lambda,\mu}(q)=0$
otherwise.
\qed
\end{theorem}

The result can also be derived as follows:
In~\cite{BL}, it is shown that the powers of the radical of $\Sym_n$ for the
internal product coincide with the homogeneous component of degree $n$ of the
lower central series of $\Sym$ for the external product:
\begin{equation}
{\mathcal R}^{*j} =\gamma^j(\Sym)
\end{equation}
where $\gamma^j(\Sym)$ is the ideal of $\NCSF$ generated by the commutators
$[\Sym,\gamma^{j-1}(\Sym)]$.

Since $\DD^{(0)}$ is a free associative algebra over a sequence of primitive
elements $(\zeta_k)_{k\ge 2}$ with the same internal product as in $\Sym$, the
argument of \cite{BL} can be reproduced \emph{verbatim}, and we see that 
\begin{equation}
({\rm rad\,}\DD^{(0)})^{*j} =\gamma^j(\DD^{(0)})= 
{\mathcal R}^{*j}\cap \DD^{(0)}\,. 
\end{equation}

This shows that, in $\DD^{(0)}$ as well as in $\Sym$, for $\lambda$ finer than
$\mu$, $e_\mu* e_I$ is nonzero modulo ${\rm rad }^{*2}$ iff $\mu$ is obtained
from $\lambda$ by summing two distinct parts.
And more generally, $e_\mu* e_I$ is in ${\rm rad}^{*k}$ and nonzero modulo
${\rm rad}^{*k+1}$ iff $\ell(\lambda)-\ell(\mu)=k$.
%
\subsection{Tables of the $q$-Cartan invariants of $\DD_n^{(0)}$}

The labels for row and columns of the $q$-Cartan matrices, namely partitions
with no part one, are in reverse lexicographic order.
In the following matrices, the zero entries are represented by dots to enhance
readability.

For $n\leq4$, the $q$-Cartan matrix of $\DD_n^{0}$ is the identity matrix of
size $d_n^{(0)}$.
The first non-trivial example arises for $n=5$.
The $q$-Cartan matrices of $\DD_5^{(0)}$, $\DD_6^{(0)}$, and $\DD_7^{(0)}$
are respectively 
{
\begin{equation}
\label{qCar-dn0-567}
\left(
\begin{array}{cc}
1 & q \\
. & 1 \\
\end{array}
\right)
\qquad
\left(
\begin{array}{cccc}
1 & q & . & . \\
. & 1 & . & . \\
. & . & 1 & . \\
. & . & . & 1 \\
\end{array}
\right)
\qquad
\left(
\begin{array}{ccccccccccccccc}
1 & q & q &q^2 \\
. & 1 & . & q  \\
. & . & 1 & .  \\
. & . & . & 1  \\
\end{array}
\right)
\end{equation}
}
and those of $\DD_8^{(0)}$ and $\DD_9^{(0)}$ are
{
\begin{equation}
\label{qCar-dn0-89}
\left(
\begin{array}{ccccccc}
1 & q & q & . &q^2&q^2& . \\
. & 1 & . & . & q & . & . \\
. & . & 1 & . & . & q & . \\
. & . & . & 1 & . & . & . \\
. & . & . & . & 1 & . & . \\
. & . & . & . & . & 1 & . \\
. & . & . & . & . & . & 1 \\
\end{array}
\right)
\qquad
\left(
\begin{array}{cccccccc}
1 & q & q & q &q^2&2q^2& . &q^3\\
. & 1 & . & . & q & q  & . &q^2\\
. & . & 1 & . & . & q  & . & . \\
. & . & . & 1 & . & q  & . & . \\
. & . & . & . & 1 & .  & . & q \\
. & . & . & . & . & 1  & . & . \\
. & . & . & . & . & .  & 1 & . \\
. & . & . & . & . & .  & . & 1 \\
\end{array}
\right)
\end{equation}
}

On these matrices, one can read the quiver of $\DD^{(0)}_n$. Note
that it is a subquiver of the quiver of $\Sym_n$ (see \cite{Sal}), since
one cannot create parts $1$ by merging parts of non-unitary partitions. 

\subsection{The (large) derangement algebra $\DD=\DD^{(\infty)}$}

Pursuing the analogy with the peak algebra, let us define
\begin{equation}
\DD := \bigoplus_{k\ge 0}S_k\DD^{(0)}\,.
\end{equation}
Note already that $\DD$ is not a subalgebra of $\NCSF$. It is only a
sub-coalgebra. Moreover, we have:

\begin{theorem}
Each homogeneous component $\DD_n$ of $\DD$ is stable under $*$. It is a
unital algebra, since it contains $S_n$, the neutral element of $*$ in
$\NCSF$.
\end{theorem}

We shall prove a slightly more general result, interpolating between
$\DD^{(0)}$ and $\DD$.

\subsection{A filtration of $\DD$}

Define $\DD_n^{(k)}$ by
\begin{equation}
\label{filt}
\DD_n^{(k)} := \bigoplus_{j=0}^k S_j \DD^{(0)}_{n-j}.
\end{equation}
For $k=0$, this is $\DD_n^{(0)}$ and for
$k\geq n$, one recovers $\DD_n$.

The following alternative definition of the filtration
will be useful in the sequel:
\begin{lemma}
\label{altfilt}
\begin{equation}
\DD_n^{(k)} = \bigoplus_{j=0}^{\min(n,k)} S_1^j \DD_{n-j}^{(0)}.
\end{equation}
\end{lemma}

\Proof
Expanding $\sigma_1=e^{S_1}\sigma_1^\sharp$, we see that
\begin{equation}
S_k\equiv \frac{S_1^k}{k!} \pmod{\bigoplus_{j<k} S_j \DD_{n-j}^{(0)} }.
\end{equation}
\qed

\subsection{Dimensions of the $\DD_n^{(k)}$}

From~(\ref{filt}), we see that the dimension $d_n^{(k)}$ of $\DD_n^{(k)}$ is 
\begin{equation}
d_n^{(k)} = \sum_{i\leq k} d_{n-i}^{(0)}.
\end{equation}
For $k=\infty$, these are the usual Fibonacci numbers.

The first values are given in the following table:
\begin{equation}
\begin{array}{c|c|c|c|c|c|c|c|c|c|c|}
  n            & 0 & 1 & 2 & 3 & 4 & 5 &  6 & 7 & 8 & 9 \\
\hline
d_n^{(0)}      & 1 & 0 & 1 & 1 & 2 & 3 &  5 & 8 & 13& 21\\
\hline
d_n^{(1)}      & 1 & 1 & 1 & 2 & 3 & 5 &  8 & 13& 21& 34\\
\hline
d_n^{(2)}      & 1 & 1 & 2 & 2 & 4 & 6 & 10 & 16& 26& 42\\
\hline
d_n^{(3)}      & 1 & 1 & 2 & 3 & 4 & 7 & 11 & 18& 29& 47\\
\hline
d_n^{(\infty)} & 1 & 1 & 2 & 3 & 5 & 8 & 13 & 21& 34& 55\\
\hline
\end{array}
\end{equation}

\subsection{The complete picture}

\begin{theorem}
For all $n$ and $k$, each homogeneous component $\DD^{(k)}_n$ of $\DD^{(k)}$
is stable under $*$.
It is a unital algebra and the neutral element is
\begin{equation}
\label{defPnk}
P_{n}^{(k)} := \sum_{i=0}^{\min(k,n)}  \frac{S_1^i}{i!} S_{n-i}^\sharp.
\end{equation}
\end{theorem}

The theorem is a consequence of the following lemma:

\begin{lemma}
\label{lem-exp}
Let $f$ and $g$ be in $\NCSF$. Then
\begin{equation}
   \left(\frac{S_1^m}{m!}f^\sharp\right)
 * \left(\frac{S_1^n}{n!} g^\sharp\right) =
\left\{
\begin{array}{cc}
0 & \text{if $m\not=n$}, \\
\frac{S_1^n}{n!} ( f^\sharp * g^\sharp) & \text{otherwise}.
\end{array}
\right.
\end{equation}
\end{lemma}

\Proof
Replacing the right factor by its generating series, we have
\begin{equation}
\begin{split}
(S_1^m f^\sharp) * (e^{S_1} g^\sharp) 
&= \sum_{(g)} (S_1^m * (e^{S_1}g_{(1)}^\sharp))\,
              (f^\sharp * (e^{S_1}g_{(2)}^\sharp)) \\
&= \sum_{(g)} (e^{S_1}g_{(1)}^\sharp * S_1^m)\,
              (f^\sharp * (e^{S_1}g_{(2)}^\sharp)), \\
\end{split}
\end{equation}
since $S_1^m$ is central for $*$.
Now, $(e^{S_1}{S^I}^\sharp * S_1^m)=0$ if $I$ is not empty since
$S_n^\sharp * S_1^n=S_1^n * S_n^\sharp=0$ for $n\geq1$ as
$e^{-S_1} * \sigma_1^\sharp=1$.
So the whole sum reduces to
\begin{equation}
\begin{split}
 S_1^m (f^\sharp * (e^{S_1}g^\sharp)) \\
\end{split}
\end{equation}
Thus,
\begin{equation}
 f^\sharp*\left(e^{S_1}g^\sharp\right)
=f*\sigma_1^\sharp*\left(e^{S_1}g^\sharp\right)
\end{equation}
Consider now the generic case $g=\sigma_1(XA)$:
\begin{equation}
\begin{split}
\sigma_1^\sharp * \left(e^{S_1}\sigma_1(XA)^\sharp\right)
&= \left(e^{-S_1} \sigma_1\right) *\left(e^{S_1}\sigma_1(XA)^\sharp\right) \\
&= \left(e^{-S_1}*e^{S_1}\sigma_1(XA)^\sharp\right)
\left(e^{S_1}\sigma_1(XA)^\sharp\right) \\
&=
  \left(e^{S_1}*e^{-S_1}\right)
  \left(\sigma_1(XA)^\sharp*e^{-S_1}\right)
  \left(e^{S_1}\sigma_1(XA)^\sharp\right) \\
&=e^{-S_1}
  \left(e^{-S_1}*\sigma_1(XA)^\sharp\right)
  \left(e^{S_1}\sigma_1(XA)^\sharp\right) \\
&=e^{-S_1}
  e^{-S_1(XA)^\sharp}
  e^{S_1}\sigma_1(XA)^\sharp \\
&=e^{-S_1} e^{S_1}\sigma_1(XA)^\sharp = \sigma_1(XA)^\sharp \\
\end{split}
\end{equation}
where the third and fourth equalities come from the fact that
$e^{-S_1}$ is central for $*$ and idempotent.
Multiplying by $f$ on the left yields
\begin{equation}
f^\sharp * \left(e^{S_1} g^\sharp\right) = f * g^\sharp,
\end{equation}
whence the statement. 
\qed

\Proof[of the theorem]
From Lemma \ref{altfilt}, we have $P_{n}^{(k)}\in\DD^{(k)}_n$.
Moreover, since the  $S_1^k {S^I}^\sharp$ with $I$ non-unitary form a basis 
of $\DD_n^{(k)}$ adapted to the direct sum decomposition, that $P_{n}^{(k)}$
is neutral on both sides is equivalent to the already known fact that
$S_n^\sharp$ is neutral in $\DD_n^{(0)}$.
\qed

\begin{corollary}
\label{cor-morph}
The map $\phi_{m} : \ \DD^{(0)}_k\rightarrow \DD^{(m)}_{k+m}$ defined by
\begin{equation}
\phi_m(f) =  \frac{S_1^{m}}{m!}f
\end{equation}
is a (non-unital) monomorphism of algebras.
\end{corollary}

\begin{corollary}
As an \emph{algebra}, $\DD_n$ is isomorphic to the direct sum
\begin{equation}
\DD_n \sim \bigoplus_{k=0}^n \DD_{n-k}^{(0)}.
\end{equation}
\end{corollary}

\Proof
As a vector space, $\DD_n$ is the direct sum of the spaces
$V_n^{(k)}=S_1^k\DD_{n-k}^{(0)}$.
By Lemma~\ref{lem-exp}, $V_n^{(k)}*V_n^{(\ell)}=0$ for $k\not=\ell$, and by
Corollary~\ref{cor-morph}, each $V_n^{(k)}$ is a subalgebra isomorphic to
$D_{n-k}^{(0)}$.
\qed

\subsection{Representation theory of $\DD_n^{(k)}$}

We are now in a position to deduce the representation theory of the
$\DD^{(k)}_n$ from that of $\DD^{(0)}_n$.
We can extend Corollary~\ref{cor-NST} to all values of $k$, as a direct
consequence of Lemma~\ref{lem-NST}.

\begin{corollary}
Let $e_I$ be the idempotents of $\NCSF$ defined in Corollary~\ref{cor-NST}.

\noindent
(i) The $e_J$ such that $J=(1^j,I)$ with $j\le k$ and $I$ does not contain a
part $1$ form a basis of $\DD_n^{(k)}$.

\noindent
(ii) The principal idempotents of $\DD^{(k)}_n$ are the $e_\lambda$
such that $m_1(\lambda)\le k$.
\qed
\end{corollary}

We can now state the general result on the representation theory of
$\DD_n^{(k)}$:

\begin{corollary}
The irreducible representations of $\DD_n^{(k)}$ are of dimension $1$ and
parametrized by partitions of $n$ with at most $k$ parts equal to $1$.
\qed
\end{corollary}

Note that the principal idempotents of $\DD_n$ are also a complete
set of minimal orthogonal idempotents of $\Sym_n$. 

\subsection{$q$-Cartan matrix of $\DD_n^{(k)}$}

We order the labels of the rows and columns of the $q$-Cartan matrices, namely
partitions with at most $k$ ones, first by their number of ones and then in
reverse lexicographic order.
So, for example, with $n=5$ and $k=3$, the order is
\begin{equation}
[ 5,\ 32,\ 41,\ 221,\ 311,\ 2111 ].
\end{equation}
With this convention, the $q$-Cartan matrix of $\DD_n^{(k)}$ is the
block-diagonal matrix obtained by putting on the diagonal the $q$-Cartan
matrices of $\DD_{n-i}^{(0)}$ for $i\leq\min(k,n)$.
Indeed, from the previous results, one easily sees that
\begin{lemma}
\begin{equation}
\DD_n^{(k)} = P_n^{(k)} * \NCSF_n * P_n^{(k)}.
\end{equation}
\end{lemma}

Since each $P_n^{(k)}$ is a sum of orthogonal idempotents (the
$S_1^iS_{n-i}^\sharp/i!$), this proves that the $q$-Cartan matrix of
$\DD_n^{(\infty)}$ is deduced from the $q$-Cartan matrix of $\NCSF_n$ by
putting to $0$ the entries whose row and column do not have the same number of
ones.

\subsection{Block projectors}

We have seen that each space $S_1^j \DD_{n-j}^{(0)}$ has as unit a part of
$P_n^{(k)}$, namely
\begin{equation}
\label{dnk}
D_{n,k} = \frac{S_1^k}{k!} S_{n-k}^\sharp
\end{equation}
Their double generating series is
\begin{equation}
D(t,u;A):= \sum_{n\ge 0}\sum_{k=0}^n t^{n-k}u^k D_{n,k}(A)
= e^{(u-t)S_1(A)}\sigma_t(A),
\end{equation}
that is, the noncommutative analog of the generating function of $d_{n,k}$,
the number of permutations in $\SG_n$ with exactly $k$ fixed points:
\begin{equation}
D(t,u)=\sum_{n\ge 0}\frac{1}{n!}\sum_{k=0}^n d_{n,k} {u^k}t^{n-k} 
      =\frac{e^{u-t}}{1-t}.
\end{equation}

\subsection{$q$-dimension polynomials}

If one sums up the entries of the $q$-Cartan matrix of $\DD_n^{(\infty)}$, one
gets the following polynomials in $q$, refining the Fibonacci numbers:
\begin{equation}
 1,\ 2,\ 3,\ 5,\ q+7,\ 2\,q+11,\ q^2+5\,q+15,\
  3\,q^2+9\,q+22, \dots 
\end{equation}
better represented in the following triangle:
\begin{equation}
\begin{array}{cccccccccc}
1   \\
2   \\
3   \\
5   \\
7   & 1   \\ 
11  & 2   \\
15  & 5   & 1   \\ 
22  & 9   & 3   \\ 
30  & 17  & 7   & 1   \\ 
42  & 28  & 16  & 3   \\ 
56  & 47  & 31  & 9   & 1   \\
77  & 73  & 58  & 21  & 4   \\
101 & 114 & 102 & 47  & 12  & 1  \\
135 & 170 & 175 & 94  & 32  & 4  \\
176 & 253 & 286 & 183 & 74  & 14 & 1  \\
231 & 365 & 461 & 333 & 162 & 40 & 5  \\
\end{array}
\end{equation}
The first column (constant terms of the polynomials) 
corresponds to the size of the matrix, hence the number of partitions of $n$.
The second column is sequence~A000097 of~\cite{Slo}:
from the characterization of the quiver of $\NCSF_n$, we also have that the second column
gives to the number of ways of selecting two different parts different from
$1$, in all partitions of $n$.
Finally, it is also equal to the number of ways of selecting two different
parts, in all partitions of $n-2$, hence justifying that the number of
arrows in the quiver of $\DD_n^{(\infty)}$ is equal to the number of arrows in
the quiver of $\NCSF_{n-2}$.



\section{The $(1-\E)$ transform in $\WQSym^*$}

\subsection{Word quasi-symmetric functions}

A word $u$ over $\N^*$ is said to be \emph{packed} if the set of letters
occuring in $u$ is an interval of $\N^*$ containing $1$.
The algebra $\WQSym(A)$ (Word Quasi-Symmetric functions) is defined as the
subalgebra of $\K\<A\>$ based on \emph{packed words} and spanned by the
elements
\begin{equation}
\M_u(A) := \sum_{\pack(w)=u} w,
\end{equation}
where $\pack(w)$ is the \emph{packed word} of $w$, that is, the word obtained
by replacing all occurrences of the $k$-th smallest letter of $w$ by $k$.
For example,
\begin{equation}
\pack(871883319) = 431442215.
\end{equation}
This is the invariant algebra of the quasi-symmetrizing action of $\SG(A)$ on
$\K\<A\>$ \cite{NCSF7}.
Packed words can be identified with set compositions in an obvious way, and
geometrically, they can be interpreted as facets of the permutohedron:
a packed word $w = w_1\dots w_n$ with largest entry $\ell$ can be identified
with the set composition $[P_1, \dots, P_\ell]$ where
$P_j = \{ i\leq n \ |\ w_i = j\}$.
For example, $431442215$ corresponds to
$[\{3,8\},\{6,7\},\{2\},\{1,4,5\},\{9\}]$.

\medskip
Let $\NW_u=\M_u^*$ be the dual basis of $(\M_u)$.
It is known that $\WQSym$ is a self-dual Hopf algebra~\cite{Hiv,NT06} and that
on the graded dual $\WQSym^*$, an internal product $*$ may be defined by
\begin{equation}
\label{intWQ}
\NW_u * \NW_v = \NW_{\pack(u,v)},
\end{equation}
where the packing of biwords is defined with respect to the lexicographic
order on biletters, so that, for example,
\begin{equation}
\pack\left(\gf{42412253}{53154323}\right)
= 62513274.
\end{equation}

This product is induced from the internal product of parking
functions~\cite{NTpark,NT1,NTp2} and allows one to identify the homogeneous
components $\WQSym_n$ with the (opposite) Solomon-Tits algebras, in the sense
of~\cite{Patras}.

The (opposite) Solomon descent algebra, realized as $\Sym_n$, is embedded in
the (opposite) Solomon-Tits algebra realized as $\WQSym^*_n$ by
\begin{equation}
\label{S2NW}
S^I = \sum_{\ev(u)=I} \NW_u,
\end{equation}
where $\ev(u)$ is the evaluation of $u$.

From now on, we shall denote $\WQSym^*$ by $\WW$.

\subsection{Idempotents of $\WW$}

\subsubsection{The semi-simple quotient}

It is known that the radical of $\WW_n$ is spanned by the differences
\begin{equation}
\NW_u - \NW_v,
\end{equation}
where $v=\sigma(u)$ for some permutation $\sigma$ of the support of $u$,
${\rm supp}(u)=\{i\ |\ |u|_i\not=0\}$.
This is easily seen: Equation~(\ref{intWQ}) implies that the $\NW_u-\NW_v$ are
nilpotent of order $2$ if $v=\sigma(u)$, and that their span is an ideal
${\mathcal R}_n$. Moreover, any product of $n$ such factors
$\NW_{u_i}-\NW_{v_i}$ vanishes since the product of two such factors is either
strictly finer than $u_i$ or zero, so that ${\mathcal R}_n$ is nilpotent.
The quotient $W_n/{\mathcal R}_n$ is semi-simple. It can be
identified with the (commutative) algebra of set partitions with $\wedge$ (the
inf for the refinement order on set partitions) as product.
Indeed, packed words encode set compositions, $u=u_1\dots u_n$ corresponding
to the set composition of $[n]$ in which $i$ belongs to the block $u_i$,
\emph{e.g.},
\begin{equation}
u = 21231 \Longleftrightarrow [\{2,5\},\{1,3\},\{4\}].
\end{equation}
and the left action of permutations amounts to permuting the blocks,
\emph{e.g.}, with $\sigma=231$,
\begin{equation}
\sigma(21231)=32312 \Longleftrightarrow [\{4\},\{2,5\},\{1,3\}].
\end{equation}

Hence, the idempotents of a complete family of $\WW_n$ are parametrized by set
partitions of $[n]$.

\subsubsection{The idempotents of Saliola}

In~\cite{Sal}, Saliola has given a general recipe for constructing such
complete sets.
Given a packed word $u$, denote by $\Pi(u)$ the set partition obtained by
forgetting the order among the blocks of the set compositions encoded by $u$.

For each set partition $\pi$ of $[n]$, choose a linear combination
\begin{equation}
l_\pi=\sum_{\Pi(u)=\pi} c_u \NW_u,
\end{equation}
where the coefficient $c_u$ depends only on the evaluation $\ev(u)$ of $u$, and
\begin{equation}
\sum_{\Pi(u)=\pi} c_u = 1.
\end{equation}
Start with the initial condition
\begin{equation}
e_{\{1\},\{2\},\dots,\{n\}}
= \frac{1}{n!} \sum_{\sigma\in\SG_n} \NW_\sigma,
\end{equation}
hence equal to $S_1^n/n!$ in~$\NCSF$,
and define by the induction
\begin{equation}
\label{recFS}
e_\pi = l_\pi * (\NW_{1^n} - \sum_{\pi'>\pi} e_{\pi'})
\end{equation}
where $\NW_{1^n}=S_n$ is the identity of $*$, and $\geq$ is the refinement
order. Then, the $e_\pi$ form a complete set of orthogonal idempotents of
$\WW_n$.

\subsubsection{A non-recursive construction}

Families of Saliola idempotents can be computed for all $\WW_n$
simultaneously, in a non-recursive way from families of idempotents of the
descent algebra $\NCSF_n$ constructed by the method developed in~\cite{NCSF2}.
Recall that the starting point of this construction is a sequence of Lie
idempotents $\gamma_n\in\NCSF_n$, that is, an arbitrary sequence of primitive
elements whose commutative image in $Sym$ is $p_n/n$.

Then, if we decompose the identity $S_n$ of $\NCSF_n$ as
\begin{equation}
S_n = \sum_{I\vDash n} c_I \gamma^I,
\end{equation}
the elements
\begin{equation}
e_\lambda := \sum_{I\downarrow n} c_I \gamma^I,
\end{equation}
with $\lambda$ a partition of $n$ form a complete family of orthogonal
idempotents of $\NCSF_n$.

Let us fix such a family, and define for each set partition $\pi$ of $[n]$
\begin{equation}
l_\pi = \sum_{\Pi(u)=\pi} c_{\ev(u)} \NW_u.
\end{equation}
These elements satisfy Saliola's conditions: obviously, $c_{\ev(u)}$ depends
only on $\ev(u)$, and
\begin{equation}
\sum_{\Pi(u)=\pi} c_{\ev(u)} = 
\prod_i m_i(\lambda)! \sum_{I\downarrow\lambda} c_I
\end{equation}
is the coefficient of $p_\lambda / z_\lambda$ in the commutative image of
$S_n$, which is $h_n$, so it is $1$.
Hence, the sequence $(\gamma_n)$ determines idempotents $e_\pi$ of the
$\WW_n$ by the recursion~(\ref{recFS}).
But we can also compute these directly as follows:

\begin{theorem}
The idempotents $e_\pi$ are given by the internal products
\begin{equation}
e_\pi = l_\pi * e_\lambda.
\end{equation}
\end{theorem}

\Proof
Let $\tilde{e}_\pi = l_\pi * e_\lambda$.
For $\pi=\{\{1\},\dots,\{n\}\}$,
we have $l_\pi = S_1^n/n!$, $e_\lambda= S_1^n/n!$, so that
$\tilde{e}_\pi = S_1^n/n!$.

Let $l_\lambda = \sum_{\Lambda(\pi)=\lambda} l_\pi$, where $\Lambda(\pi)$ is
the integer partition recording the block lengths of $\pi$.
We have
\begin{equation}
l_\lambda = \sum_{I\downarrow \lambda} c_I S^I
\equiv e_\lambda \mod\bigoplus_{l(I)>l(\lambda)} \K \gamma^I,
\end{equation}
so that $e_\lambda = l_\lambda * e_\lambda$ in $\NCSF$.

We want to show that $\tilde{e}_\pi = e_\pi$.
For that, recall from~\cite{Sal} that, if $\Pi(u)\not\leq\pi$,
$\NW_u*e_\pi=0$,
so that
$l_\pi*e_{\pi'}=0$ if $\pi\not>\pi'$, and
$\Lambda(\pi')\not=\Lambda(\pi)$.
This implies in particular that
$e_\lambda = \sum_{\Lambda(\pi)=\lambda} e_\pi$.
Indeed, this is true for $\lambda=1^n$, and, by induction,
$e_\pi = l_\pi* ( \NW_{1^n} - \sum_{l(\pi')>l(\pi)} e_{\pi'})$,
since $\pi'>\pi$ implies $l(\pi')>l(\pi)$.
Hence
\begin{equation}
e_\pi = l_\pi* ( \NW_{1^n} - \sum_{l(\lambda')>l(\lambda)} e_{\lambda'})
= l_\pi* \sum_{l(\lambda')\leq l(\lambda)} e_{\lambda'}.
\end{equation}
Summing over $\pi$, we get
\begin{equation}
\sum_{\Lambda(\pi)=\lambda} e_\pi
= l_\lambda * \sum_{l(\lambda')\leq l(\lambda)} e_{\lambda'}
= e_\lambda.
\end{equation}

Now,
\begin{equation}
\begin{split}
e_\pi 
&= l_\pi * (\NW_{1^n} - \sum_{\pi'>\pi} e_{\pi'})
= l_\pi * \sum_{\pi'\not>\pi} e_{\pi'}
\\
&= l_\pi * (e_\lambda + \sum_{\pi'\not>\pi ; \Lambda(\pi')\not=\Lambda(\pi)}
            e_{\pi'})
 = l_\pi * e_\lambda = \tilde{e}_\pi.
\end{split}
\end{equation}
\qed

\subsection{The $(1-\E)$-transform in $\WW$}

The embedding~(\ref{S2NW}) of $\Sym$ in $\WW$ can be defined on the
generators as
\begin{equation}
S_n\longmapsto \NW_{1^n}\,.
\end{equation}
It is clearly a bialgebra morphism.
The element
\begin{equation}
\sigma_1^\sharp = e^{-\NW_1}\sum_{n\ge 0}\NW_{1^n}
\end{equation}
is well-defined in $\WW$, and so is the $\sharp$-transform
\begin{equation}
F^\sharp := F * \sigma_1^\sharp\,.
\end{equation}

\subsection{Bases of $\WW^\sharp$}

Let us say that a packed word $u$ is {\em non-unitary} (and \emph{unitary}
otherwise) if no letter occurs exactly once in $u$. These words correspond to
set compositions without singletons.

\begin{proposition}
The $\NW_u^\sharp$ with u non-unitary form a basis of $\WW^\sharp$.
\end{proposition}

\Proof
Let us say that $v$ is finer than $u$ (and write $v>u$) if the set composition
encoded by $v$ is finer than the set composition encoding $v$.
Then,
\begin{equation}
\NW_u * \sigma_1^\sharp = \NW_u + \sum_{v} {c_{uv}\NW_v},
\end{equation}
where $v>u$ or $v$ is unitary. Hence, the $\NW_v^\sharp$ with $u$ non-unitary
are linearly independent.
\qed

\subsection{Algebraic structure of $\WW^\sharp$}

Let ${\mathcal J}$ be the two-sided ideal of $\WW$ generated by the $\NW_u$
such that $u$ has at least a letter occuring exactly once.
The product rule (\ref{intWQ}) shows that ${\mathcal J}$ is an ideal for the
internal product as well. Hence, the projection
\begin{equation}
\pi:\ \WW\longrightarrow \WW/{\mathcal J}
\end{equation}
is a morphism for $*$. Its restriction to $\WW^\sharp$ is then an isomorphism,
and clearly,
\begin{equation}
\pi(\sigma_1^\sharp)=\sigma_1\,.
\end{equation}
Since $\sigma_1$ is neutral in $\WW/{\mathcal J}$, we have:
\begin{proposition}
$\sigma_1^\sharp$ is neutral in $\WW^\sharp$.\qed
\end{proposition}

\noindent
Note that this proof would apply to $\Sym$ as well.
To summarize,

\begin{proposition}
$\WW^\sharp$ is isomorphic to $\WW/{\mathcal J}$ as a Hopf algebra, and each
$\WW_n^\sharp$ is $*$-isomorphic to $\WW_n/{\mathcal J}_n$, with
$D_n=S_n^\sharp=\NW_{1^n}^\sharp$ as neutral element.
\qed
\end{proposition}

\subsection{Representation theory of $\WW^\sharp$}

We can now apply Lemma \ref{BAe} with 
$A=\WW_n$,
$B=\WW_n^\sharp$ and $e=\S_n^\sharp$.

Thus
\begin{equation}
\rad^k \WW_n^\sharp = D_n * (\rad^k \WW_n) * D_n.
\end{equation}

The irreducible representations of $\WW_n$, which are one-dimensional
are parametrized by set partitions of $[n]$.

The $q$-Cartan matrices and quiver of $\WW_n$ have been determined
in~\cite{Sc}:

\begin{equation}
\label{cabq}
c_{\alpha,\beta}(q) = c_{\alpha,\beta} q^{l(\alpha)-l(\beta)},
\end{equation}
where $l(\pi)$ is the number of blocks of a set partition $\pi$, and the
Cartan invariant $c_{\alpha,\beta}$ is $0$ if $\alpha$ is not finer than
$\beta$, and otherwise
\begin{equation}
c_{\alpha,\beta} = \prod_i (m_i-1)!,
\end{equation}
where for each block $B_i$ of $\beta$, $m_i$ is the number of blocks of
$\alpha$ into which $B_i$ has been split.

For example, with $\alpha = 12|3|4|56|7$ and $\beta=1234|567$,
we have: $c_{\alpha,\beta} = (3-1)! (2-1)! = 2$, $l(\alpha)=5$,
$l(\beta)=2$, so that
$c_{\alpha,\beta}(q) = 2q^3$.

\begin{theorem}
The $q$-Cartan matrix of $\WW_n^\sharp$ is the restriction to rows and columns
indexed by non-unitary set partitions of $[n]$ of~(\ref{cabq}).
In particular, the vertices of the quiver of $\WW_n^\sharp$ are the
non-unitary set partitions, and there is an arrow $\alpha\to\beta$ whenever
$\beta$ is obtained from $\alpha$ by merging two blocks.
\end{theorem}


\subsection{Analogue of $\DD$ in $\WW$}

Let ${\mathcal V}_n^{(0)}=\WW_n^\sharp$ and
\begin{equation}
{\mathcal V}_n^{(k)} = \bigoplus_{k=0}^n  D_{n,k} * \WW_n * D_{n,k}\,.
\end{equation}

Then, as in the case of $\NCSF$, each ${\mathcal V}_n^{(k)}$ is a unital
subalgebra of $\WW_n^\sharp$.

%
%
\section{The $(1-\E)$-transform in $\FQSym$}

\subsection{Definition}

Recall that $\FQSym$ is based on permutations, that in the mutually dual
bases $\F_\sigma=\G_{\sigma^{-1}}$, the internal product is defined by
\begin{equation}
\F_\sigma*\F_\tau=\F_{\sigma\tau}\quad\text{or equivalently}\quad
\G_\sigma*\G_\tau=\G_{\tau\sigma}\,,
\end{equation}
and that $\Sym$ is embedded into $\FQSym$ by $S_n=\G_{12\dots n}$.
The transformation can therefore be defined by
\begin{equation}
\F_\sigma^\sharp := \F_\sigma*\sigma_1^\sharp\,.
\end{equation}
Since the splitting formula remains valid in $\FQSym$ when the right factor
of the internal product is in $\Sym$ \cite{NCSF7}, this is again a Hopf
algebra morphism.

As we shall see below, in $\FQSym$, the idempotent $D_n=S_n^\sharp$ as well
as the other $D_{n,k}$ defined in~(\ref{dnk}) admit an interesting
interpretation.

\subsection{The Tsetlin library (uniform case)}

The $(1-\E)$-transform in $\FQSym$ is related to a classical problem in
probability theory known as the \emph{Tsetlin library} (see {\it e.g.,}
\cite{BHR}).
This is a Markov chain on $\SG_n$, defined by a shelf of $n$ books, which are
randomly picked by users and them put back at the left of the shelf after use.
In the uniform case (when all books are picked with the same probability), the
determination of the stationary distribution amounts to the diagonalisation of
the linear operator on $\C\SG_n$
\begin{equation}
t_n(f)=f\tau_n
\end{equation}
where
\begin{equation}
\tau_n = 12\ldots n+ 2134\dots n + 3124\dots n +\cdots + n12\dots n-1
 \in \C\SG_n\,.
\end{equation}
This problem is also an ingredient of the proof of Hivert's conjecture by
Garsia and Wallach~\cite{GW}. It can be solved in many different ways.
The following one is quite natural in the context on Noncommutative
Symmetric Functions.

We start with the observation that $\tau_n$ is in the descent algebra
$\Sigma_n$. Indeed, $\tau_n=D_{\subseteq \{1\}}$ (the sum of permutations
having at most a descent at the first position), so that its representation
as a noncommutative symmetric function is $S^{1,n-1}$, a rather
well-understood element.

From this remark, we obtain immediately the eigenvalues of $t_n$. Indeed,
according to Proposition 3.12 of \cite{NCSF2}, these are the scalar products
$\langle h_{n-1}h_1,p_\lambda\rangle$ of ordinary symmetric functions.
Clearly, the scalar product evaluates to $m_1(\lambda)$, so that the spectrum
is $0,1,2,\ldots,n-2,n$.

Let us now construct the spectral projectors. To this aim, we shall need to
evaluate some polynomials in $t_n$.
Let us set $T_n=S^{1,n-1}$ and consider the generating function
\begin{equation}
T=\sum_{n\ge 0}T_n = S_1\sigma_1\,.
\end{equation}
Since the internal products $T_n*T_m$ are (by definition) $0$ for
$m\not=n$, we have
\begin{equation}
\sum_n T_n^{*r} = T^{*r}\,
\end{equation}
and using iteratively the splitting formula (\cite{NCSF2}, Prop. 2.1)
\begin{equation}
T*T^{*(r-1)} = (S_1\sigma_1)*T^{*(r-1)}=\mu[(S_1\otimes\sigma_1) *_2
\Delta(T^{*(r-1)})],
\end{equation}
we get the expression
\begin{equation}
T^{*r} =  B_r(S_1)\sigma_1,
\end{equation}
where $B_r(x)$ are the Bell polynomials (this is the obvious noncommutative
analogue of the classical formula for the Kronecker powers of the
representation of $\SG_n$ by permutation matrices).

Using the fact that the coefficients of $B_n$ are the Stirling numbers of the
second kind $S(n,k)$, we obtain in $\Sym$
\begin{equation}
T*(T-1)*(T-2)*\cdots * (T-k+1)= S_1^k\sigma_1\,,
\end{equation}
and in particular, in degree $n$,
\begin{equation}
T_n*(T_n-1)*(T_n-2)*\cdots * (T_n-n+1) = S_1^n
\end{equation}
and as well
\begin{equation}
T_n*(T_n-1)*(T_n-2)*\cdots * (T_n-n+2) =S_1^n\,,
\end{equation}
and since it is plain that $S_1^n*T_n=nS_1^n$, so that $S_1^n*(T_n-n)=0$,
the minimum polynomial of $T_n$ is
\begin{equation}
P_n(x)=(x-n)\prod_{k=0}^{n-2}(x-k)\,.
\end{equation}
This shows that $T_n$ is semisimple, and allows an easy construction
of the spectral projectors.

Let us start with the kernel. The projector is given by $f\mapsto f*D_n$
where
\begin{equation}
D_n=\frac{(T_n-1)*(T_n-2)*\cdots*(T_n-n+2)*(T_n-n)}{(-1)(-2)\cdots(2-n)(-n)}
\end{equation}
but since $T_n-n+1$ is invertible, one can take as well
\begin{equation}
D_n=\frac{(-1)^n}{n!}(T_n-1)_{*n}
\end{equation}
where $(x)_n=x(x-1)\cdots (x-n+1)$, and the star means evaluation with the
internal product. This is a better choice, since we have now a simple
generating series for all these projectors,
\begin{equation}
\sum_{n\ge 0}D_nx^n=e^{-xS_1}\cdot\sigma_x= \sigma_x((1-\E)A)\,.
\end{equation}
Indeed, we have $(x-1)_n=(x)_n-n(x-1)_{n-1}$, so that
\begin{equation}
\frac{(-1)^n}{n!}(T-1)_{*n} = \frac{(-S_1)^n}{n!}\sigma_1-
\frac{(-1)^{n-1}}{(n-1)!}(T-1)_{*\
n-1}=\sum_{k=0}^n\frac{(-S_1)^k}{k!}\sigma_1\,.
\end{equation}

The same reasoning shows that the projectors $D_{n,k}$ on the eigenspaces of
$k$ are given by the generating series
\begin{equation}
D(t,u)=\sum_{n\ge 0}\sum_{k=0}^nt^{n-k}u^k D_{n,k} = e^{(u-t)S_1}\sigma_1 \,,
\end{equation}
which is $\sigma_1(tA-(t-u)\E A)$, so that these elements coincide with those
defined by~(\ref{dnk}).

\subsection{Characters of the associated modules}

The Frobenius characteristic of the left ideal of $\C\SG_n$ generated by
the idempotents $\delta_{n,k}$ corresponding to $D_{n,k}$ via the
identification $\sigma\leftrightarrow \F_\sigma$ can now be calculated
as follows (compare \cite[Cor. 4.2]{Sc}).

Since $\delta_{n,k}$ is an idempotent, its characteristic $\ch(\delta_{n,k})$
coincides with its cycle index $Z(\delta_{n,k})$. By the Gessel-Reutenauer
formula \cite{GR}, the coefficient of $p_\mu$ in  $Z(\delta_{n,k})$ is equal
to $\< \underline{D_{n,k}}, L_\mu\>$, where $\underline{F}$ means the
commutative image of the noncommutative symmetric function $F$,
and for $\nu=1^{n_1}2^{n_2}\ldots$,
\begin{equation}
L_\nu=h_{n_1}[\ell_1]h_{n_2}[\ell_2]\cdots\,,\quad
\ell_n=\frac1n\sum_{d|n}\mu(d)p_{n/d}^d\,,
\end{equation}
$\mu$ denoting here the Moebius function.
Hence, the generating function of all the cycle indexes is
\begin{equation}
Z_Y(t,u)=\<D(t,u;X),  L(X,Y)  \>_X
=\<D(t,u;X),L(Y,X)\>_X
\end{equation}
by the symmetry formula \cite{ST}
\begin{equation}
L(X,Y)=L(Y,X)=\sum_\mu L_\mu(X)p_\mu(Y) = \prod_{n\ge
1}\sigma_{p_n(X)}[\ell_n(Y)]\,.
\end{equation}
Plugging this last expression into the scalar product and dualizing, we
obtain
\begin{equation}
  \prod_{n\ge 1}\sigma_{p_n(t+(u-t)\E)}[\ell_n(Y)]
  = \frac{\sigma_1((u-t)Y)}{1-tp_1(Y)}\,.
\end{equation}
In particular, specializing at $Y=\E$ gives that the dimension
of $\C\SG_n\delta_{n,k}$ is $d_{n,k}$, the number of permutations in
$\SG_n$ with exactly $k$ fixed points.

It is also easy to obtain the expansion of $Z_Y(t,u)$ as a combination
of the $L_\mu(Y)$. Indeed, writing $D(t,u;x)=\sigma_1((u-t)\E X +tX)$,
we have $\<D(t,u),p_\mu\>=u^{m_1}\prod_{i\ge 2}t^{m_i}$, where $m_i$
is the multiplicity of the part $i$ in $\mu$. Hence,
\begin{equation}
  Z(\delta_{n,k})=\sum_{m_1(\mu)=k}L_\mu
\end{equation}
We note that this is the quasi-symmetric generating function of the
permutations with exactly $k$ fixed points.
Note that $Z(D(t)$ is the commutative image of the generating series
of desarrangements (\ref{desar}).

\subsection{$q$-derangement numbers}

From the above considerations, one can easily derive a (known) closed formula
for the $q$-derangement numbers (compare \cite[Theorem 4.5]{Sc})
\begin{equation}
d_n(q):=\sum_{\sigma\in D_n}q^{\maj \sigma}\,.
\end{equation}
Indeed,
\begin{equation}
\begin{split}
d_n(q)&=\left\langle \sum_{\sigma\in D_n}\F_\sigma, 
        \sum_{\tau\in\SG_n}q^{\maj(\tau)}\G_\tau\right\rangle_{\FQSym}\\
&=\sum_{\sigma\in D_n}\sum_{|I|=n}q^{\maj(I)}\<F_{C(\sigma)},R_I\>
=\left\langle
\sum_{m_1(\mu)=0}L_\mu, K_n(q)\right\rangle\,. 
\end{split} 
\end{equation}
Hence,
\begin{equation}
\begin{split}
\sum_{n\ge 0} d_n(q)\frac{x^n}{(q)_n}
&=\left\langle \sigma_1\left[\sum_{n\ge 2}x^nl_n\right] ,
        \sigma_1\left(\frac{X}{1-q}\right)\right\rangle\\
&=\left\langle \frac{\lambda_{-x}}{1-xp_1},
              \sigma_1\left(\frac{X}{1\!-\!q}\right)
 \right\rangle
=\lambda_{-x}\left(\frac{1}{1\!-\!q}\right)
 \left(1-\frac{x}{1\!-\!q}\right)^{-1}\\
&=\frac{1-q}{1-x-q}\prod_{n\ge 0}(1-xq^n)
\end{split}
\end{equation}
so that finally~\cite{Wa}
\begin{equation}
d_n(q)=[n]!\sum_{k=0}^n\frac{(-1)^k}{[k]!}q^{\binom{k}{2}}\,.
\end{equation}

\subsection{Characters from Lie idempotents}

The expression
\begin{equation}
\label{charact}
{\rm ch}(\delta_{n,k})=\sum_{m_1(\lambda)=k}L_\lambda\,,
\end{equation}
is also a consequence of the following type of expressions
\begin{equation}
\label{decD}
D_{n,k}=\sum_{m_1(\lambda)=k}E_\lambda(\pi)
\end{equation}
in the notation of \cite[Theorem 3.16]{NCSF2}, for some sequence $\pi_n$ of
Lie idempotents in descent algebras. Indeed, \cite[Theorem 3.21]{NCSF2},
implies then that the character is given by (\ref{charact}). We have already
seen one such expression with $\pi_n=\zeta_n$, the Zassenhaus idempotents.
We can also write one involving the Hausdorff series. Writing as usual
\begin{equation}
\Phi =\sum_{n\ge 1}\phi_n=\log \sigma_1
\end{equation}
(the Solomon idempotents), we have
\begin{equation}
\sum_{n\ge 0}D_n=  e^{-\phi_1}  e^\Phi=e^{H(-\phi_1,\Phi)}
\end{equation}
where $H$ is the Hausdorff series. Taking $\pi_1=\phi_1$ and
$\pi_n=H_n(-\phi_1,\Phi)$ for $n\ge 2$, we obtain a sequence
of Lie idempotents (see, \emph{e.g.}, \cite{NCSF2}, Theorem 3.1),
from which it is easy to build a decomposition of the identity
\begin{equation}
\sigma_1 = e^{\pi_1}\exp\left\{\sum_{n\ge 2}\pi_n\right\}\,,
\end{equation}
and more explicitely,
\begin{equation}
\label{decId}
S_n=\sum_{r+s=n}\frac{1}{r!s!}\sum_{\ell(J)=r,|J|=n-s,1\not\in J}\pi^{1^s,J}.
\end{equation}
This gives in particular the decomposition (\ref{decD}),
with, for a partition $\lambda$ such that $m_1(\lambda)=0$,
\begin{equation}\label{projpi}
E_\lambda(\pi) =\frac{1}{\ell(\lambda)!}\sum_{I\downarrow \lambda}\pi^I\,,
\end{equation}
where $I\downarrow \lambda$ means that the nondecreasing rearrangement
of the composition $I$ is the partition $\lambda$.

\subsection{Eigenbases of $t_n$}

From Proposition 7.4 of \cite{NCSF3}, we know the image a projector
of the type (\ref{projpi}).
It is formed of weighted symmetrizations of Lie elements.
With the above $\pi_n$, the distribution is uniform, so that the kernel
consists in ordinary symmetrized products of Lie elements. Concretely, a basis
of $\Ker t_n$ in $\C\SG_n$ is for example
\begin{equation}
(\gamma_1\theta_{\lambda_1},\gamma_2\theta_{\lambda_2},\cdots,
\gamma_r\theta_{\lambda_r})
\end{equation}
where $(a,b,c)=abc+acb+bac+bca+cab+cba$, and so on (symmetrized products),
the $\gamma_k$ are the minimal representatives of the cycles of a
derangement, $\theta_n=[[\cdots[1,2],3],\cdots n]$ is a Dynkin element,
and $\lambda$ runs over partitions without part 1.

For example, a basis of $\Ker t_4$ of dimension $d_4=9$ is given by the
elements $[[[1,a],b]c]$ with $abc$ running over permutations of $234$, for
$L_4$, and by the three symmetrized products $([1,2],[3,4])$, $([1,3],[2,4])$,
and $([1,4],[2,3])$ for $L_{22}$.

Bases of the other eigenspaces are obtained by the same process, using
weighted symmetrizations as indicated in~\cite{NCSF3}.
Indeed, Equation~(\ref{decId}) shows that a basis of the eigenspace with
eigenvalue $s$ is given by
\begin{equation}
(\gamma_1\theta_{\lambda_1},\gamma_2\theta_{\lambda_2},\cdots,
\gamma_r\theta_{\lambda_r})\cdot (j_1\shuffle j_2\shuffle\cdots\shuffle j_s)
\end{equation}
where $\gamma_k$ are the minimal representatives of the cycles of length at
least $2$ a permutation of cycle type $(\lambda,1^s)$ having $s$ fixed points
$j_1,j_2,\ldots,j_s$.

To continue with $n=4$, a basis of the $1$ eigenspace is $[[i,j],k]\cdot l$
($i<j,k$, $ijkl$ a permutation of $1234$), dimension $8$, and a basis of the
$2$-eigenspace is given by $$[i,j]\cdot(kl+lk)\,\ i<j,\ k<l$$ where $ijkl$ is
a permutation of $1234$ (dimension $6$).
Finally, the $4$-eigenspace is one dimensional and generated by the full
symmetrizer.

Using the $\zeta_n$ instead of the $\pi_n$, we can replace symmetrized
products by ordinary products of homogeneous Lie polynomials taken
in nondecreassing order of the degrees.

The idempotents $\delta_{n,k}$ have  been first studied by M.
Schocker~\cite{Sc} (apparently unaware of previous works on the subject and of
their relation with the Tsetlin library).

\subsection{A basis of $\FQSym^\sharp$}

We have seen that the $(1-\E)$-transform is a bialgebra morphism in $\FQSym$.
Hence, its image $\FQSym^\sharp$ is a Hopf subalgebra.
The $\SG_n$-module $\Delta_{n,k}$ can be identified with
\begin{equation}
\FQSym_n^\sharp = \FQSym_n * D_{n,k},
\end{equation}
so that
\begin{equation}
\dim \FQSym_n^\sharp = d_n.
\end{equation}
It is therefore desirable to find a basis of $\FQSym^\sharp$ labeled by
derangements, or some other set of permutations naturally in bijection with
these. As we shall see, the natural transformation involved here is simply a
version of Foata's first fundamental transformation~\cite{Loth}. 

Let $\gamma_n$ be the cycle
\begin{equation}
\gamma_n := n\,1\,2\dots n-1,
\end{equation}
so that
\begin{equation}
\S^{\gamma_n} = T_n=S_1 S_{n-1} = R_n + R_{1,n-1} =
\sum_{\sigma\in 1\shuffle 23\dots n} \F_\sigma.
\end{equation}
Since the $\sharp$-transform is a morphism for the product of $\FQSym$,
\begin{equation}
\S^{\gamma_n} * \sigma_1^\sharp = S_1^\sharp S_{n-1}^\sharp =0,
\end{equation}
and, for any permutation $\sigma\in\SG_n$,
\begin{equation}
 (\F_\sigma * \S^{\gamma_n})^\sharp
= \F_\sigma * \S^{\gamma_n} * \sigma^\sharp
= 0.
\end{equation}
Recall that $i$ is a \emph{left-right minimum} of $\sigma$ if
\begin{equation}
\sigma_j>\sigma_i \text{\ for all $j<i$}.
\end{equation}
Let $X_n$ be the set of permutations of $\SG_n$ such that $\sigma\cdot 0$ does
not have two consecutive left-right minima (that is, $\sigma$ does not end by
$1$ and does not have two consecutive LR-minima),
and let $Y_n = \SG_n \backslash X_n$.

\begin{lemma}
\label{FY}
For $\sigma\in Y_n$, write
\begin{equation}
\sigma\cdot 0 = \cdots \sigma_i \sigma_{i+1} \cdots
\end{equation}
where $i$ and $i+1$ is the first pair of consecutive LR-minima, and let
\begin{equation}
\sigma' = \sigma_i \cdot \sigma_1 \cdots \hat{\sigma_i} \cdots \sigma_n
\end{equation}
be the permutation obtained by moving $\sigma_i$ at the first position,
leaving the remaining letters unchanged, and removing the zero in the end.
Then
\begin{equation}
\F_{\sigma'} * \S^{\gamma_n} = \F_\sigma + \sum_{\tau\in T} \F_\tau,
\end{equation}
where the permutations of $T$ are lexicographically smaller than $\sigma$.
\end{lemma}

\Proof
The expression
\begin{equation}
\F_{\sigma'} * \S^{\gamma_n} =
\sum_{\tau \in \sigma_i\shuffle \sigma_1 \dots \hat{\sigma_i} \dots \sigma_n}
\F_\tau
\end{equation}
contains $\F_\sigma$ and, since $\sigma_i$ is a LR-minimum, $\sigma$ is the
maximal element of the previous sum.
\qed

\medskip
For a permutation $\sigma$ with LR-minima $i_1,\dots,i_p$, let
\begin{equation}
\phi(\sigma) = (\sigma_1\dots\sigma_{i_1-1})
(\sigma_{i_1}\dots\sigma_{i_2-1})
\dots
(\sigma_{i_p}\dots\sigma_{n}),
\end{equation}
where each parenthesis represents a cycle.
For example, with $\sigma = 62781453$,
\begin{equation}
\phi(\sigma) = (6)(278)(1453) = 47153682.
\end{equation}
This is Foata's first fundamental transformation (up to reversing the order on
the integers), hence a bijection.
Clearly, $\phi(\sigma)$ has fixed points whenever $\sigma\in Y_n$, so
$\phi$ induces a bijection between $X_n$ and derangements of $\SG_n$.

From Lemma~\ref{FY}, we see that the elements
\begin{equation}
(\F_\sigma^\sharp)_{\sigma\in X_n}
\end{equation}
span $\FQSym_n^\sharp$. Since $|X_n|= d_n = \dim \FQSym_n^\sharp$, we have
finally
\begin{theorem}
The $\F_\sigma^\sharp$ for $\sigma\in X_n$ form a basis of $\FQSym^\sharp$.
\qed
\end{theorem}

The sets $X_n$ have an interesting structure.

\begin{theorem}
The set $X_n$ is an ideal of the left weak order on $\SG_n$.
Its maximal elements are the left-shifted concatenations
\begin{equation}
w_I := w_{i_1} \btr \dots \btr w_{i_r},
\end{equation}
where $w_i:=1\,i\,i\!-\!1,\dots,2$, composition $I$ has no part $1$, and
$\alpha\btr\beta = \alpha[\ell]\cdot \beta$ if $\beta\in\SG_\ell$.
\end{theorem}

\Proof
To show that $X_n$ is an ideal, we will prove that if $s_i$ denotes the
elementary transposition $(i,i+1)$, then $\sigma\in Y_n$ and
$\inv(s_i\sigma)=\inv(\sigma)+1$, implies $s_i\sigma\in Y_n$.

If $\sigma_k=r > s=\sigma_{k+1}$ are consecutive LR-minima of $\sigma$,
they will remain so for $s_i\sigma$, unless $i=r-1$, $r$, $s-1$, or $s$.
Since $s_i\sigma$ has one inversion more than $\sigma$, we can exclude the
case $i=r-1$: $r$ being a LR-minimum, $r-1$ cannot be to the left of $r$ in
$\sigma$. We can also exclude $i=s-1$ for the same reason.
If $i=r$, then $r$ is exchanged with $r+1$, which has to be to its right in
$\sigma$, so that again $\sigma_k$ and $\sigma_{k+1}$ are consecutive
LR-minima in $s_i\sigma$.
The same reasoning applies with $i=s$.
Hence $Y_n$ is a coideal, and consequently $X_n$ is an ideal.

Now, the elements $w_I$ are clearly in $X_n$ when $I$ has no part $1$, and any
exchange of consecutive values creating an inversion in such a $w_I$ would
create a pair of consecutive LR-minima. So these $w_I$ are maximal elements of
the ideal $X_n$.

Conversely, consider $\sigma\in X_n$ maximal.
Then consider the suffix $s$ of $\sigma$ beginning with $1$. The
maximality condition of $\sigma$ implies that if $t$ belongs to that suffix,
then $t-1$ also belongs to it. So this prefix is a permutation of an
$\SG_{|s|}$, then should be $1|s|\dots 2$. The same now works by induction on
the permutation $\tau$ defined by $\sigma = \tau\btr s$.
\qed

\medskip
For example, with $n=5$, we get the following three maximal elements of $X_n$:
\begin{equation}
15432,\ 35412,\ 45132.
\end{equation}

\medskip
The same proof can be adapted to the case of permutations that are images by
$\phi$ of permutations with at most $k$ fixed points.
Let $X_n^{(k)}$ be the image by $\phi$ of permutations with at most $k$ fixed
points. Then $X_n^{(k)}$ is the set of permutations with at most $k$
consecutive LR-minima.

\begin{theorem}
The set $X_n^{(k)}$ is a ideal of the left weak order on $\SG_n$.
Its maximal elements are the $w_I$ where $I$ runs over compositions with
\begin{itemize}
\item either $k-1$ ones and the remaining parts equal to $2$,
\item or exactly $k$ ones.
\end{itemize}
\end{theorem}

\Proof
The fact that $X_n^{(k)}$ is an ideal comes from the same idea as before: all
permutations greater than a given permutation $\sigma$ for the left weak order
have LR-minima at the same position.

By the same argument as in the previous theorem, the maximal elements must be
some $w_I$, where $I$ has at most $k$ ones. Now, it is clear that
\begin{equation}
w_I < w_J,
\end{equation}
in the left weak order iff $I$ can be obtained from $J$ by gluing parts equal
to $1$ with their next part. So the compositions described in the statement
are definitely maximal elements. And since all compositions with at most $k$
ones can be obtained from these ones by the gluing process, this ends the
proof.
\qed

\bigskip
Here is a table of the number of maximal elements of $X_n^{(k)}$

\begin{equation}
\begin{array}{|c||c|c|c|c|c|c|c|c|c|}
\hline
n\backslash k & 0 & 1 & 2 & 3 & 4 & 5 & 6 & 7 & 8 \\
\hline
\hline
1             & 0 & 1 &   &   &   &   &   &   &   \\
\hline
2             & 1 & 1 & 1 &   &   &   &   &   &   \\
\hline
3             & 1 & 2 & 2 & 1 &   &   &   &   &   \\
\hline
4             & 2 & 3 & 3 & 3 & 1 &   &   &   &   \\
\hline
5             & 3 & 5 & 6 & 4 & 4 & 1 &   &   &   \\
\hline
6             & 5 & 9 & 9 &10 & 5 & 5 & 1 &   &   \\
\hline
7             & 8 &15 &19 &14 &15 & 6 & 6 & 1 &   \\
\hline
8             &13 &27 &31 &34 &20 &21 & 7 & 7 & 1 \\
\hline
\end{array}
\end{equation}

Note that the first column is obviously given by Fibonacci numbers since these
indeed count the number of compositions of $n$ in parts at least $2$.
The other columns are not known to~\cite{Slo} and neither is the sequence of
row sums.

But there exists a simple formula giving the number of maximal elements,
coming directly from their characterization:
\begin{equation}
m_{n,k} := \binom{(n+k-1)/2}{k-1} +
           \sum_{\ell=0}^{\lfloor[\frac{n-k}{2}\rfloor]}
            \binom{\ell+k}{k} \binom{n-k-\ell-1}{\ell-1},
\end{equation}
with the convention that a binomial coefficient with entries not in the
natural numbers is zero.

\subsection{Other bases of $\FQSym^\sharp$}


\begin{conjecture}
Let $<'$ be the order on permutations defined by
\begin{equation}
\sigma <' \tau \Longleftrightarrow \phi(\sigma) <_{\rm lex} \phi(\tau).
\end{equation}
Then the matrix of ${\S^\sigma}^\sharp$ of the $\S$ basis is triangular.
Moreover, the diagonal values are $1$ for the elements of $X_n$ and $0$ for
$Y_n$.
\end{conjecture}

\bigskip
For example, here are the matrices for $n=2$, $3$, and $4$ (Figure 1) where the zero
entries have been represented by dots to enhance readability.

The  permutations are ordered as follows:
\begin{equation}
[21,\ 12],\qquad\qquad [321,\ 312,\ 231,\ 123,\ 132,\ 213].
\end{equation}
\begin{equation}
\begin{split}
[& 4321,\ 4312,\ 4231,\ 4123,\ 4132,\ 4213,\ 3421,\ 3412,\\
 & 2341,\ 1234,\ 1243,\ 2314,\ 2431,\ 1423,\ 3241,\ 2134,\\
 & 3142,\ 1324,\ 1432,\ 2413,\ 2143,\ 3214,\ 1342,\ 3124].
\end{split}
\end{equation}

\begin{equation}
\left(\begin{array}{cc}
. & -1/2 \\
. & 1    \\
\end{array}
\right)
\qquad
\left(\begin{array}{cccccc}
. & . & .  & 1/3 & -1/3 & 2/3 \\
. & . & .  & -1  &  .   & -1  \\
. & . & .  &  .  &  .   &  .  \\
. & . & .  &  1  &  .   &  1  \\
. & . & .  &  .  &  1   & -1  \\
. & . & .  &  .  &  .   &  .  \\
\end{array}
\right)
\end{equation}

\begin{figure}[ht]
{\tiny
\rotateleft{
$
\left(\begin{array}{cccccccccccccccccccccccccccccc}
.& .& .& .& .& .& .&  1/4& .& -1/8& 1/4& -3/8& .&  1/8& .& -1/4&  3/8&
-1/4& -1/4&  3/8& 1/2& -3/4&  1/8& -3/8  \\
.& .& .& .& .& .& .& -1/2& .&  1/2&  .&   1/2& .& -1/2& .&   1&    .&
1/2&   .&  -1/2&  .&    1&    .&    1   \\
.& .& .& .& .& .& .&   .&  .&   .&   .&    .&  .&   .&  .&   .&    .&
.&    .&    .&   .&    .&  -1/2&   .   \\
.& .& .& .& .& .& .&   .&  .&  -1&   .&   -1&  .&   .&  .&  -1&    .&
-1/2&   .&    .&   .&   -1&    .&   -1   \\
.& .& .& .& .& .& .&   .&  .&   .&   -1&   1&  .&   .&  .&   .&   -1&
1/2&   .&    .&   -1&   1&    .&    .   \\
.& .& .& .& .& .& .&   .&  .&   .&   .&    .&  .&   .&  .&   .&    .&
-1/2&   .&    .&   .&    .&    .&    .   \\
.& .& .& .& .& .& .& -1/2& .&   .&   .&    .&  .&   .&  .&   .&    .&
.&    .&    .&   .&    .&    .&    .   \\
.& .& .& .& .& .& .&   1&  .&   .&   .&    .&  .&   .&  .&   .&    .&
.&    .&    .&   .&    .&    .&    .   \\
.& .& .& .& .& .& .&   .&  .&   .&   .&    .&  .&   .&  .&   .&   1/2&
.&    .&    .&   .&    .&    .&    .   \\
.& .& .& .& .& .& .&   .&  .&   1&   .&    1&  .&   .&  .&   1&    .&
.&    .&    .&   .&    1&    .&   1/2  \\
.& .& .& .& .& .& .&   .&  .&   .&   1&   -1&  .&   .&  .&   .&    .&
.&    .&    .&   1&   -1&    .&  -1/2  \\
.& .& .& .& .& .& .&   .&  .&   .&   .&    .&  .&   .&  .&   .&    .&
.&    .&    .&   .&    .&    .&   1/2  \\
.& .& .& .& .& .& .&   .&  .&   .&   .&    .&  .&   .&  .&   .&  -1/2&
.&    .&    .&   .&    .&    .&    .   \\
.& .& .& .& .& .& .&   .&  .&   .&   .&    .&  .&   1&  .&  -1&    .&
.&    .&    1&   .&   -1&    .&  -1/2  \\
.& .& .& .& .& .& .&   .&  .&   .&   .&    .&  .&   .&  .&   .&  -1/2&
.&    .&    .&   .&    .&    .&    .   \\
.& .& .& .& .& .& .&   .&  .&   .&   .&    .&  .&   .&  .&   .&    .&
.&    .&    .&   .&    .&    .&  -1/2  \\
.& .& .& .& .& .& .&   .&  .&   .&   .&    .&  .&   .&  .&   .&    1&
.&    .&    .&   .&    .&    .&    .   \\
.& .& .& .& .& .& .&   .&  .&   .&   .&    .&  .&   .&  .&   .&    .&
1&    .&    .&   .&    .&    .&    .   \\
.& .& .& .& .& .& .&   .&  .&   .&   .&    .&  .&   .&  .&   .&    .&
.&    1&   -1&   -1&   1&    .&   1/2  \\
.& .& .& .& .& .& .&   .&  .&   .&   .&    .&  .&   .&  .&   .&    .&
.&    .&    .&   .&    .&    .&  -1/2  \\
.& .& .& .& .& .& .&   .&  .&   .&   .&    .&  .&   .&  .&   .&    .&
.&    .&    .&   .&    .&    .&   1/2  \\
.& .& .& .& .& .& .&   .&  .&   .&   .&    .&  .&   .&  .&   .&    .&
.&    .&    .&   .&    .&    .&  -1/2  \\
.& .& .& .& .& .& .&   .&  .&   .&   .&    .&  .&   .&  .&   .&    .&
.&    .&    .&   .&    .&    1&    .   \\
.& .& .& .& .& .& .&   .&  .&   .&   .&    .&  .&   .&  .&   .&    .&
.&    .&    .&   .&    .&    .&    1   \\
\end{array}
\right)
$
}
}
\caption{The  matrix of $\sharp$ in the $\S$ basis of $\FQSym$}
\end{figure}

\section{Other combinatorial Hopf algebras}

\subsection{The algebras $\PQSym$ and $\CQSym$}

There is an internal product on $\PQSym$ extending that of
$\WQSym^*$~\cite{NTpark}.  The $\sharp$-transform is defined in $\PQSym$
(it contains $\NCSF$ as a subalgebra), but
$\PQSym*\WQSym^* \subseteq \WQSym^*$, so that $\PQSym^\sharp = \WW^\sharp$,
and we get nothing new.

\medskip
Similarly, the Catalan algebra $\CQSym$~\cite{NTp2}, we have
\begin{equation}
\CQSym * \NCSF \subseteq \NCSF,
\end{equation}
so that
\begin{equation}
\CQSym^\sharp = \NCSF^\sharp.
\end{equation}

\subsection{The algebra of planar binary trees $\PBT$}

The Loday-Ronco algebra of planar binary tree is not stable by the
$\sharp$-transform. Since $\PBT$ is the subalgebra of $\FQSym$ generated by
the $S^{\sigma}$ where $\sigma$ avoids the pattern $132$ (see~\cite{HNT}), we
have, for example :
\begin{equation}
{S^{213}}^\sharp = S^{123} - S^{132} - S^{312} + \frac{2}{3} S^{321}
\not\in \PBT.
\end{equation}
However, $\PBT^\sharp$ is a well-defined Hopf subalgebra of $\FQSym$.

\begin{conjecture}
The algebra $\PBT^\sharp$ is free over the set $\P_T^\sharp$, where $T$ runs
over trees with at least two nodes, and such that the right subtree of the root
is empty.
\end{conjecture}

In particular, the conjecture implies that the dimension of the homogeneous
components $\PBT_n^\sharp$ are given by the Fine numbers~\cite{DS,Slo},
sequence A000957:
\begin{equation}
1,\,0\, 1\, 2\, 6\, 18\, 57\, 186,\, 622 \dots
\end{equation}


\newpage
\footnotesize
\begin{thebibliography}{abc}
%
\bibitem{ANT}{\sc M. Aguiar, J.-C. Novelli} and {\sc J.-Y. Thibon},
{\it Unital versions of the higher order peak algebras},
arXiv:0810.4634, FPSAC'09 (Linz).
%
\bibitem{BHR}{\sc P. Bidigare, P.  Hanlon} and {\sc D.  Rockmore},
{\it  A combinatorial description of the spectrum of the Tsetlin library and
its generalization to hyperplane arrangements},
Duke Math. J. {\bf 99} (1999), 135--174.
%
\bibitem{BHT}{\sc N. Bergeron, F. Hivert} and {\sc J.-Y. Thibon},
{\it The peak algebra and the Hecke-Clifford algebras at $q=0$ },
J. Combinatorial Theory A {\bf 117} (2004), 1--19.
%
\bibitem{BL}{\sc D.Blessenohl} and {\sc H. Laue}
{\it The module structure of Solomon's descent algebra},
J. Aust. Math. Soc. {\bf 72} (2002), no. 3, 317--333.
%
\bibitem{D1}{\sc J. D\'esarm\'enien}, 
{\it Une autre interpr\'etation du nombre de d\'erangements},
S\'emin. Lothar. Comb. {\bf 8} (1984) 11--16. 
%
\bibitem{DW1}{\sc J. D\'esarm\'enien} and {\sc M. Wachs},
{\it Descentes sur les d\'erangements et mots circulaires}, 
S\'emin. Lothar. Comb. {\bf 19} (1988), 13--21.
%
\bibitem{DW2}{\sc J. D\'esarm\'enien} and {\sc M. Wachs},
{\it Descent classes of permutations with a given number of fixed points}, 
J. Combin. Theory Ser. A {\bf 64} (1993) 311--328.
%
\bibitem{DS}{\sc E. Deutsch} and {\sc L. Shapiro},
{\it A survey of the Fine numbers},
Discrete Math., {\bf 241} (2001), 241--265. 
%
\bibitem{NCSF6}{\sc G. Duchamp, F. Hivert}, and {\sc J.-Y. Thibon},
{\it Noncommutative symmetric functions VI: free quasi-symmetric functions and
related algebras},
Internat. J. Alg. Comput. {\bf 12} (2002), 671--717.
%
\bibitem{NCSF7}{\sc G. Duchamp, F. Hivert, J.-C. Novelli} and {\sc J.-Y.
Thibon},
{\it Noncommutative Symmetric Functions VII: Free Quasi-Symmetric Functions
Revisited}, arXiv:0809.4479.
%
\bibitem{NCSF3}{\sc G. Duchamp, A. Klyachko, D. Krob} and {\sc J.-Y. Thibon},
{\it Noncommutative symmetric functions III: Deformations of Cauchy and
convolution algebras},
Discrete Mathematics and Theoretical Computer Science {\bf 1} (1997), 159--216.
%
\bibitem{NCSF1}{\sc I.M. Gelfand, D. Krob, A. Lascoux, B. Leclerc,
V.~S. Retakh}, and {\sc J.-Y. Thibon},
{\it Noncommutative symmetric functions},
Adv. in Math. {\bf 112} (1995), 218--348.
%
\bibitem{HNT}{\sc F. Hivert, J.-C. Novelli}, and {\sc J.-Y. Thibon}.
{\it The algebra of binary search trees},
Theoret. Comput. Sci. {\bf 339} (2005), 129--165.
%
\bibitem{Loth}{\sc N. Lothaire},
{\it Combinatorics on words},
Cambridge University Press (1997).
%
\bibitem{GaR}{\sc A. M. Garsia} and {\sc C. Reutenauer},
{\it A decomposition of Solomon's descent algebra},
Advances in Math. {\bf 77} (1989), 189--262.
%
\bibitem{GW}{\sc A. M. Garsia} and {\sc N. Wallach},
{\it $r$-$QSym$ is free over $Sym$},
Jour. Comb. Theo A, {\bf 114} (2007), 704--732.
%
\bibitem{GR}{\sc I. Gessel} and {\sc C. Reutenauer},
{\it Counting permutations with given cycle structure and descent set},  
J. Comb. Theory A {\bf 64} (1993), 189-215.
%
\bibitem{Hiv} {\sc F. Hivert},
{\it Combinatoire des fonctions quasi-sym\'etriques},
Th\`ese de Doctorat, Marne-La-Vall\'ee, 1999.
%
\bibitem{NCSF2}{\sc D. Krob, B. Leclerc} and {\sc J.-Y. Thibon},
{\it Noncommutative symmetric functions II{\,}: Transformations of alphabets},
Int. J. of Alg. and Comput. {\bf 7} (1997), 181--264.
%
\bibitem{Las}{\sc A. Lascoux},
{\it Symmetric functions and combinatorial operators on polynomials},
CBMS Regional Conference Series in Mathematics {\bf 99},
American Math. Soc., Providence, RI, 2003; xii+268 pp.
%
\bibitem{LS}{\sc A. Lascoux} and {\sc M.~P. Sch\"utzenberger},
{ \it Formulaire raisonn\'e de fonctions sym\'etriques},
Publ. Math. Univ. Paris 7, Paris, 1985.
%
\bibitem{Mcd}{\sc I.G. Macdonald},
{\it Symmetric functions and Hall polynomials},
2nd ed., Oxford University Press, 1995.
%
\bibitem{MR}{\sc C. Malvenuto} and {\sc C. Reutenauer},
{\it Duality between quasi-symmetric functions and the Solomon descent
algebra}, J. Algebra {\bf 177} (1995), 967--982.
%
\bibitem{NST}{\sc J.-C. Novelli, F. Saliola}, and {\sc J.-Y. Thibon},
{\it Representation theory of the higher order peak algebras},
preprint math.CO/0906.5236.
%
\bibitem{NT1}{\sc J.-C. Novelli} and {\sc J.-Y. Thibon},
{\it A Hopf algebra of parking functions},
FPSAC'04, Vancouver, 2004.
%
\bibitem{NTpark}{\sc J.-C. Novelli} and {\sc J.-Y. Thibon},
{\it Parking functions and descent algebras},
Annals of Comb., {\bf 11} (2007), 59--68.
%
\bibitem{NT06}{\sc J.-C. Novelli} and {\sc J.-Y. Thibon},
{\it Polynomial realizations of some trialgebras},
FPSAC'06. Also preprint ArXiv:math.CO/0605061.
%
\bibitem{NTp2}{\sc J.-C. Novelli} and {\sc J.-Y. Thibon},
{\it Hopf algebras and dendriform structures arising from parking functions},
Fund. Math., 193 (2007), 189--241.
%
\bibitem{NT3}{\sc J.-C Novelli} and {\sc J.-Y. Thibon},
{\it Noncommutative symmetric functions and Lagrange inversion},
Adv. Appl. Math. {\bf 40} (2008), 8--35.
%
\bibitem{Patras}{\sc F. Patras} and {\sc M. Schocker},
{\it Twisted Descent Algebras and the Solomon-Tits Algebra},
Adv. in Math. {\bf 199} (2006), 151--184.
%
\bibitem{Sal}{\sc F. Saliola},
{\it On the quiver of the descent algebra},
J. Algebra, {\bf 320}, (2008), pp. 3866-3894.
%
\bibitem{Sc}{\sc M. Schocker},
{\it Idempotents for derangement numbers},
Discrete Math. {\bf 269} (2003), 239--248.
%
\bibitem{ST}{\sc T. Scharf} and {\sc J.-Y. Thibon},
{\it A Hopf algebra approach to inner plethysm},
Adv Math. {\bf 104} (1994), 30--58.
%
\bibitem{Slo}{\sc N. J. A. Sloane},
{\it The On-Line Encyclopedia of Integer Sequences} (electronic),
http://www.research.att.com/${}^\sim$njas/sequences/
%
\bibitem{Wa}{\sc M. Wachs}, {\it On $q$-derangement numbers},
Proc. Amer. Math. Soc. {\bf 106} (1989), 273--278.
\end{thebibliography}
\normalsize
\end{document}